
\def\LaTeX{%
  \let\Begin\begin
  \let\End\end
  \let\salta\relax
  \let\finqui\relax
  \let\futuro\relax}

\def\UK{\def\our{our}\let\sz s}
\def\USA{\def\our{or}\let\sz z}

\UK



\LaTeX

\USA


\salta

\documentclass[twoside,12pt]{article}
\setlength{\textheight}{24cm}
\setlength{\textwidth}{16cm}
\setlength{\oddsidemargin}{2mm}
\setlength{\evensidemargin}{2mm}
\setlength{\topmargin}{-15mm}
\parskip2mm


\usepackage[usenames,dvipsnames]{color}
\usepackage{amsmath}
\usepackage{amsthm}
\usepackage{amssymb}
\usepackage[mathcal]{euscript}

\usepackage{cite}
%
%


\definecolor{viola}{rgb}{0.3,0,0.7}
\definecolor{ciclamino}{rgb}{0.5,0,0.5}

\def\gianni #1{{\color{viola}#1}}
\def\pier #1{{\color{red}#1}}
\def\juerg #1{{\color{Green}#1}}

\def\pier #1{#1}
\def\juerg #1{#1}
\def\gianni #1{#1}




\bibliographystyle{plain}


%

\finqui

\def\Beq{\Begin{equation}}
\def\Eeq{\End{equation}}
\def\Bsist{\Begin{eqnarray}}
\def\Esist{\End{eqnarray}}

\def\Bthm{\Begin{theorem}}
\def\Ethm{\End{theorem}}
\def\Blem{\Begin{lemma}}
\def\Elem{\End{lemma}}

\def\Brem{\Begin{remark}\rm}
\def\Erem{\End{remark}}

\def\Bcenter{\Begin{center}}
\def\Ecenter{\End{center}}
\let\non\nonumber




\def\step #1 \par{\medskip\noindent{\bf #1.}\quad}


\def\Lip{Lip\-schitz}
\def\Holder{H\"older}

\def\aand{\quad\hbox{and}\quad}
\def\lsc{lower semicontinuous}

\def\lhs{left-hand side}
\def\rhs{right-hand side}
\def\sfw{straightforward}



\def\multibold #1{\def\arg{#1}%
  \ifx\arg\pto \let\next\relax
  \else
  \def\next{\expandafter
    \def\csname #1#1#1\endcsname{{\bf #1}}%
    \multibold}%
  \fi \next}

\def\pto{.}

\def\multical #1{\def\arg{#1}%
  \ifx\arg\pto \let\next\relax
  \else
  \def\next{\expandafter
    \def\csname cal#1\endcsname{{\cal #1}}%
    \multical}%
  \fi \next}


\def\multimathop #1 {\def\arg{#1}%
  \ifx\arg\pto \let\next\relax
  \else
  \def\next{\expandafter
    \def\csname #1\endcsname{\mathop{\rm #1}\nolimits}%
    \multimathop}%
  \fi \next}

\multibold
qwertyuiopasdfghjklzxcvbnmQWERTYUIOPASDFGHJKLZXCVBNM.

\multical
QWERTYUIOPASDFGHJKLZXCVBNM.

\multimathop
diag dist div dom mean meas sign supp .

\def\Span{\mathop{\rm span}\nolimits}


\def\accorpa #1#2{\eqref{#1}--\eqref{#2}}
\def\Accorpa #1#2 #3 {\gdef #1{\eqref{#2}--\eqref{#3}}%
  \wlog{}\wlog{\string #1 -> #2 - #3}\wlog{}}


\def\separa{\noalign{\allowbreak}}

\def\neto{\mathrel{{\scriptscriptstyle\nearrow}}}
\def\seto{\mathrel{{\scriptscriptstyle\searrow}}}

\def\somma #1#2#3{\sum_{#1=#2}^{#3}}
\def\tonde #1{\left(#1\right)}

\def\graffe #1{\mathopen\{#1\mathclose\}}

\def\<#1>{\mathopen\langle #1\mathclose\rangle}
\def\norma #1{\mathopen \| #1\mathclose \|}

\def\[#1]{\mathopen\langle\!\langle #1\mathclose\rangle\!\rangle}

\def\iot {\int_0^t}

\def\intQt{\int_{Q_t}}
\def\intQ{\int_Q}
\def\iO{\int_\Omega}
\def\iG{\int_\Gamma}
\def\intS{\int_\Sigma}
\def\intSt{\int_{\Sigma_t}}

\def\dt{\partial_t}
\def\dn{\partial_\nu}

\def\cpto{\,\cdot\,}

\def\checkmmode #1{\relax\ifmmode\hbox{#1}\else{#1}\fi}
\def\aeO{\checkmmode{a.e.\ in~$\Omega$}}
\def\aeQ{\checkmmode{a.e.\ in~$Q$}}
\def\aeG{\checkmmode{a.e.\ on~$\Gamma$}}
\def\aeS{\checkmmode{a.e.\ on~$\Sigma$}}
\def\aet{\checkmmode{a.e.\ in~$(0,T)$}}

\def\aat{\checkmmode{for a.a.~$t\in(0,T)$}}


\def\erre{{\mathbb{R}}}




\def\genspazio #1#2#3#4#5{#1^{#2}(#5,#4;#3)}
\def\spazio #1#2#3{\genspazio {#1}{#2}{#3}T0}

\def\L {\spazio L}
\def\H {\spazio H}
\def\W {\spazio W}


\def\Lx #1{L^{#1}(\Omega)}
\def\Hx #1{H^{#1}(\Omega)}

\def\LxG #1{L^{#1}(\Gamma)}
\def\HxG #1{H^{#1}(\Gamma)}

\def\Cx #1{C^{#1}(\overline\Omega)}

\def\Luno{\Lx 1}
\def\Ldue{\Lx 2}

\def\Huno{\Hx 1}
\def\Hdue{\Hx 2}
\def\Hunoz{{H^1_0(\Omega)}}
\def\HunoG{\HxG 1}
\def\HdueG{\HxG 2}

\def\LunoG{\LxG 1}
\def\LdueG{\LxG 2}



\let\badtheta\theta
\let\theta\vartheta
\let\eps\varepsilon
\let\phi\varphi

\let\TeXchi\chi                         
\newbox\chibox
\setbox0 \hbox{\mathsurround0pt $\TeXchi$}
\setbox\chibox \hbox{\raise\dp0 \box 0 }
\def\chi{\copy\chibox}


\def\QED{\hfill $\square$}


\def\CO{C_\Omega}

\def\suG{_{|\Gamma}}
\def\suS{_{|\Sigma}}

\def\VG{V_\Gamma}
\def\HG{H_\Gamma}
\def\WG{W_\Gamma}
\def\nablaG{\nabla_\Gamma}
\def\DeltaG{\Delta_\Gamma}
\def\muG{\mu_\Gamma}
\def\rhoG{\rho_\Gamma}
\def\tauO{\tau_\Omega}
\def\tauG{\tau_\Gamma}
\def\fG{f_\Gamma}
\def\vG{v_\Gamma}
\def\wG{w_\Gamma}
\def\gG{g_\Gamma}
\def\zetaG{\zeta_\Gamma}
\def\xiG{\xi_\Gamma}
\def\soluz{\juerg{((\mu,\muG),(\rho,\rhoG),(\zeta,\zetaG))}}
\def\Mu{(\mu,\muG)}
\def\Rho{(\rho,\rhoG)}
\def\Zeta{(\zeta,\zetaG)}
\def\Xi{(\xi,\xiG)}

\def\gstar{g^*}

\def\rhoz{\rho_0}
\def\rhoGz{{\rhoz}\suG}
\def\mz{m_0}

\def\mueps{\mu^\eps}
\def\rhoeps{\rho^\eps}
\def\muGeps{\mu_\Gamma^\eps}
\def\rhoGeps{\rho_\Gamma^\eps}
\def\zetaeps{\zeta^\eps}
\def\zetaGeps{\zeta_\Gamma^\eps}
\def\Mueps{(\mueps,\muGeps)}
\def\Rhoeps{(\rhoeps,\rhoGeps)}

\def\tauOeps{\tau_\Omega^\eps}
\def\tauGeps{\tau_\Gamma^\eps}

\def\Beta{\hat\beta}
\def\BetaG{\Beta_\Gamma}
\def\betaG{\beta_\Gamma}
\def\betaeps{\beta_\eps}
\def\betaGeps{\beta_{\Gamma\!,\,\eps}}
\def\Betaeps{\hat\beta_\eps}
\def\BetaGeps{\hat\beta_{\Gamma\!,\,\eps}}
\def\betaz{\beta^\circ}
\def\betaGz{\betaG^\circ}
\def\Pi{\hat\pi}
\def\PiG{\Pi_\Gamma}
\def\piG{\pi_\Gamma}

\def\mun{\mu^n}
\def\rhon{\rho^n}
\def\muGn{\mu_\Gamma^n}
\def\rhoGn{\rho_\Gamma^n}
\def\Mun{(\mun,\muGn)}
\def\Rhon{(\rhon,\rhoGn)}
\def\rhozn{\rho_0^n}

\def\ei{e^i}
\def\eGi{e_\Gamma^i}
\def\Ei{(\ei,\eGi)}
\def\ej{e^j}
\def\eGj{e_\Gamma^j}
\def\Ej{(\ej,\eGj)}
\def\muj{\mu_j}
\def\rhoj{\rho_j}
\def\mubar{\overline\mu}
\def\rhobar{\overline\rho}

\def\rhomin{\rho_*}
\def\rhomax{\rho^*}

\def\calVz{\calV_0}
\def\calHz{\calH_0}
\def\calVsz{\calV_{*0}}
\def\calVzp{\calV_0^{\,*}}
\def\calVp{\calV^{\,*}}
\def\calVn{\calV_n}

\def\normaV #1{\norma{#1}_V}
\def\normaH #1{\norma{#1}_H}
\def\normaW #1{\norma{#1}_W}
\def\normaVG #1{\norma{#1}_{\VG}}
\def\normaHG #1{\norma{#1}_{\HG}}
\def\normaWG #1{\norma{#1}_{\WG}}
\def\normaHH #1{\norma{#1}_{\calH}}
\def\normaVV #1{\norma{#1}_{\calV}}
\def\normaWW #1{\norma{#1}_{\calW}}

\def\ei{e^i}
\def\ej{e^j}
\def\eGi{e_\Gamma^i}
\def\eGj{e_\Gamma^j}

\let\hat\widehat

\Begin{document}


%
\title{On a Cahn--Hilliard system with convection\\[0.3cm] 
  and dynamic boundary conditions}
\author{}
\date{}
\maketitle
\Bcenter
\vskip-1cm
{\large\sc Pierluigi Colli$^{(1)}$}\\
{\normalsize e-mail: {\tt pierluigi.colli@unipv.it}}\\[.25cm]
{\large\sc Gianni Gilardi$^{(1)}$}\\
{\normalsize e-mail: {\tt gianni.gilardi@unipv.it}}\\[.25cm]
{\large\sc J\"urgen Sprekels$^{(2)}$}\\
{\normalsize e-mail: {\tt sprekels@wias-berlin.de}}\\[.45cm]
$^{(1)}$
{\small Dipartimento di Matematica ``F. Casorati'', Universit\`a di Pavia}\\
{\small and Research Associate at the IMATI -- C.N.R. Pavia}\\
{\small via Ferrata 5, 27100 Pavia, Italy}\\[.2cm]
$^{(2)}$
{\small Department of Mathematics}\\
{\small Humboldt-Universit\"at zu Berlin}\\
{\small Unter den Linden 6, 10099 Berlin, Germany}\\[2mm]
{\small and}\\[2mm]
{\small Weierstrass Institute for Applied Analysis and Stochastics}\\
{\small Mohrenstrasse 39, 10117 Berlin, Germany}
\Ecenter
\Begin{abstract}
This paper deals with an initial and boundary value problem for a system 
coupling equation and boundary condition both of Cahn--Hilliard type; an
additional convective term with a forced velocity field, which could act as a control 
on the system, is also present in the equation. Either regular or singular potentials 
are admitted in the bulk and on the boundary. Both the viscous and pure Cahn--Hilliard 
cases are investigated, and a number of results is proven about existence of solutions, 
uniqueness, regularity, continuous dependence, uniform boundedness of solutions, strict 
separation property. A complete approximation of the problem, based on the regularization 
of maximal monotone graphs and the use of a Faedo--Galerkin scheme, is introduced and 
rigorously discussed.  
\End{abstract}
\salta
\pagestyle{myheadings}
\newcommand\testopari{\sc Colli \ --- \ Gilardi \ --- \ Sprekels}
\newcommand\testodispari{\sc Cahn--Hilliard system with convection 
  and dynamic b.c.}
\markboth{\testopari}{\testodispari}
\finqui
%

\section{Introduction}
\label{Intro}
\setcounter{equation}{0}

This paper is concerned with the following Cahn-Hilliard system with convection:
\Beq
  \dt\rho + \nabla\rho \cdot u - \Delta\mu = 0
  \aand
  \tauO \dt\rho - \Delta\rho + f'(\rho) = \mu
  \quad \hbox{in $Q:=\Omega\times(0,T)$},
  \label{Isystem}
\Eeq
where $\Omega$ denotes a bounded three-dimensional domain and $T>0 $ is a fixed final time. 
The unknowns are~$\rho$, the order parameter, and~$\mu$, the chemical potential;
$f'$~stands for the derivative of a double-well potential~$f$,  
$u$ is a given velocity field and $\tauO$~is a nonnegative constant. 
According to whether $\tauO$~is positive or zero,
we speak of viscous Cahn--Hilliard or 
pure Cahn--Hilliard system, respectively. 

The equations in \eqref{Isystem} provide a 
description of the evolution phenomena related to solid-solid phase separations
with convection leaded by the term  $\nabla\rho \cdot u$, for some fixed velocity vector $u$. Let us refer to\cite{CahH, EllSh, NovCoh, bai, EllSt} for some pioneering contributions on the class of Cahn--Hilliard problems.  
In general, an evolution process goes on with diffusion;  
however, for the process of phase separation there is a structural difference since 
each phase concentrates and the so-called  spinodal decomposition occurs. A discussion on the modeling approach for phase separation, 
spinodal decomposition and mobility of atoms between cells 
can be found in \cite{FG, Gu, Podio, CGPS3, CMZ11}). 

Typical and important examples of $f$ 
are the so--called {\em classical regular potential} and the {\em logarithmic double-well potential\/}.
They are given~by
\begin{align}
  & f_{reg}(r) := \frac 14 \, (r^2-1)^2 \,,
  \quad r \in \erre, 
  \label{regpot}
  \\
  & f_{log}(r) := ((1+r)\ln (1+r)+(1-r)\ln (1-r)) - c r^2 \,,
  \quad r \in (-1,1),
  \label{logpot}
\end{align}
where $c>1$ \juerg{is such} that $f_{log}$ is nonconvex.
Another example is the following {\em double obstacle potential\/}:
\Beq
  f_{2obs}(r) := - cr^2 
  \quad \hbox{if $|r|\leq1$}
  \aand
  f_{2obs}(r) := +\infty
  \quad \hbox{if $|r|>1$},
  \label{obspot}
\Eeq
where $c>0$. \juerg{In cases} like \eqref{obspot}, one has to split $f$ into a non-differentiable convex part 
(the~indicator function of $[-1,1]$ in the present example) and a smooth perturbation.
Accordingly, one has to replace the derivative of the convex part
by the subdifferential and \juerg{interpret the second identity in \eqref{Isystem}} as a differential inclusion. In order to incorporate cases like \eqref{obspot} in our analysis, we allow $f'$ to be expressed by the sum 
$\beta + \pi$, where $\beta$ is the subdifferential of a convex and \lsc\ function
$\Beta : \erre \to [0, +\infty] $ such that $\Beta (0)=0$, and $\pi$ is the Lipschitz continuous derivative of the concave perturbation $\Pi  : \erre \to \erre$. Thus, we have that $f= \Beta + \Pi$ represents a possibly non-smooth double-well potential.

In order to set an initial-boundary value problem for \eqref{Isystem}, we have to specify initial and boundary conditions.
As \juerg{far as} the latter are concerned, the classical ones are
the homogeneous Neumann boundary conditions, namely
\Beq
  \dn\mu = 0, \quad \dn\rho = 0 
  \quad \hbox{ on $\Sigma:=\Gamma\times(0,T)$},
  \label{pier7}
\Eeq
where $\Gamma$ stands for  the smooth boundary of $\Omega$ and 
$\dn$ denotes the outward normal derivative. In the present work, 
on the contrary we tackle two dynamic boundary conditions for $\mu$ 
and~$\rho$ so to obtain a system of Cahn--Hilliard type also 
on the boundary. Namely, we complement \eqref{Isystem}~with
\Beq
  \dt\rhoG + \dn\mu - \DeltaG\muG = 0
  \aand
  \tauG \dt\rhoG + \dn\rho - \DeltaG\rhoG + \fG'(\rhoG) = \muG
  \quad \hbox{on $\Sigma$},
  \label{IdynBC}
\Eeq
where $\muG$ and $\rhoG$ are the traces of $\mu$ and~$\rho$, respectively,
$\DeltaG$~is the Laplace-Beltrami operator on the boundary,
$\tauG$~is a nonnegative constant,
and $\fG'= \betaG + \piG$ comes out from another potential~$\fG= \BetaG + \PiG$
with the same behavior as $f$, the two potentials being not completely independent
but related by a suitable growth condition. Then, it turns out that initial conditions
should be prescribed both in the bulk and on the boundary.

Therefore, by considering everything, the resulting initial and boundary value 
problem~reads
\begin{align}
  & \dt\rho + \nabla\rho \cdot u - \Delta\mu = 0
  \quad \hbox{in $Q$},
  \label{pier2}
    \\
  & \tauO \dt\rho - \Delta\rho + \beta(\rho) + \pi(\rho) \ni \mu
  \quad \hbox{in $Q$},
  \label{pier3}
  \\
  &\rhoG= \rho_{|_\Sigma}, \quad  \muG= \mu_{|_\Sigma} \aand  \dt\rhoG + \dn\mu - \DeltaG\muG = 0
  \quad \hbox{on $\Sigma$},
  \label{pier4}
  \\
  & \tauG \dt\rhoG + \dn\rho - \DeltaG\rhoG + \betaG(\rhoG) + \piG(\rhoG) \ni \muG
  \quad \hbox{on $\Sigma$},
  \label{pier5}
  \\
 & \rho(0) = \rhoz
  \quad \hbox{in $\Omega$} \aand \rhoG (0) = \rhoGz \quad \hbox{on $\Gamma$}.
    \label{pier6}
\end{align}
Up to our knowledge, in the case of a pure Cahn--Hilliard system, 
that is, with $\tauO=\tauG=0$, and without convective term ($u=0$), 
the problem \eqref{pier2}--\eqref{pier6} has been firstly formulated  
by {G}oldstein, {M}iranville and {S}chimperna~\cite{GMS11} 
and analyzed from various viewpoints in other contributions 
(see~\cite{CGM13, CMZ11, CP14, GM13}); moreover, in the case 
of general potentials, the problem has been deeply investigated in 
\cite{CF2} from the point of view of existence, uniqueness and 
regularity of the weak solution (see also \cite{FY} for an optimal 
control problem) by using an abstract approach. Here, instead, we 
face with {\it the full system} \eqref{pier2}--\eqref{pier6} by a 
complete approximation procedure, which involves not only a regularization 
of graphs but the setting of a precise Faedo--Galerkin scheme. Moreover, 
in the viscous case with both $\tauO$ and $\tauG$ positive, we can prove 
the uniform boundedness of both the chemical potential and the order 
parameter, up to the boundary, and we are even able to show the strict 
separation property in the case of logarithmic potentials like $f_{log}$ 
in \eqref{logpot}. In addition to this, we did our best to try to keep 
minimal assumptions on the velocity field $u$, concerning summability 
and time derivation (see the later~\eqref{hpu} and~\eqref{hpureg}). 
So, we think that our contribution could be a useful tool for studying 
other problems, which possibly involve other equations with coupled terms, 
and in particular for investigating optimal control problems. 

Let us now review some related literature. It turns out that some 
class of Cahn--Hilliard system, possibly including dynamic boundary conditions, 
has collected a noteworthy interest in recent years: we can quote 
\cite{CFP, Kub12, MZ, PRZ, RZ, WZ} among other contributions. In case 
of no convective term in \eqref{pier2}, and assuming the homogeneous 
boundary condition $\dn \mu =0$ (i.e., the first condition in \eqref{pier7}) 
and the  condition
\eqref{pier5} with $\tauG>0$ and $\muG$ as a given datum, the problem 
has been first addressed in \cite{GiMiSchi}: the well-posedness and 
the large time behavior of solutions have been studied for regular 
potentials $f$ and $f_\Gamma$, as well as for various singular 
potentials like the ones in \eqref{logpot} and \eqref{obspot}. 
One can see \cite{GiMiSchi, GiMiSchi2}: in these two papers the authors 
were able to overcome the difficulties due to singularities using a set 
of assumptions for $\beta, \, \pi$ and $\betaG, \, \piG $ that gives the 
role of the dominating potential to $f$ and entails some technical difficulties.  
The subsequent papers 
\cite{CGS3, CGS5, CGS4} follow a different approach 
(firstly considered in \cite{CaCo} to investigate the Allen--Cahn 
equation with dynamic boundary conditions), which consists in 
letting $f_\Gamma$ be the leading potential with respect to $f$: 
by this the analysis turns out to be simpler. The paper 
\cite{CGS3}~contains many results about existence, uniqueness 
and regularity of solutions for general potentials that include \accorpa{regpot}{logpot},
and are valid for both the viscous and pure cases, i.e., by assuming just $\tauO\geq0$.
Moreover, the optimal boundary control problems for the viscous 
and pure Cahn--Hilliard equation are discussed in \cite{CGS4} 
and \cite{CGS5}, respectively, in analogy with the corresponding 
contributions \cite{CS,CFS} for the Allen--Cahn equation. The paper 
\cite{CF1} deals with the well-posedness of the same system, but in 
which also an additional mass constraint on the boundary is imposed. 
In addition, we aim to emphasize that Cahn--Hilliard systems have been rather 
investigated from the viewpoint of optimal control. In this connection, 
we point out the contributions \cite{ZL1, ZL2} dealing with the convective 
Cahn--Hilliard equation;  the case with a nonlocal potential is studied in \cite{RoSp}. We also 
refer to \cite{CFGS1, HW1, WaNa, ZW} and quote the paper \cite{CFGS2} investigating the 
second-order optimality conditions for the state system considered 
in  \cite{CGS4}. There also exist articles addressing some discretized versions 
of general Cahn--Hilliard systems, cf. \cite{HW2, Wang}. 

The present paper is organized as follows. 
In the next two sections, we list our assumptions and notations, state our results
and give the relations between weak solutions and the above boundary value problem.
Sections~\ref{CONTDEP} is devoted to continuous dependence and uniqueness,
while the existence of a solution is shown in Section~\ref{EXISTENCE}
by taking the limit of suitable approximating problems studied in Section~\ref{APPROXIMATION}.
Finally, Section~\ref{COMPLEMENTS} is devoted to our regularity results.


\section{Statement of the problem and results}
\label{STATEMENT}
\setcounter{equation}{0}

In this section, we state precise assumptions and notations and present our results.
First of all, the \juerg{set} $\Omega\subset\erre^3$ is assumed to be bounded, connected and smooth.
As in the Introduction, $\nu$~is the outward unit normal vector field on $\Gamma:=\partial\Omega$, 
and $\dn$ and $\DeltaG$ stand for the corresponding normal derivative
and the Laplace-Beltrami operator, respectively.
Furthermore, we denote by $\nablaG$ the surface gradient
and write $|\Omega|$ and $|\Gamma|$ 
for the volume of $\Omega$ and the area of~$\Gamma$, respectively.

If $X$ is a Banach space, $\norma\cpto_X$ denotes both its norm and the norm of~$X^3$. 
Moreover, $X^*$~is the dual space of $X$, and $\<\cpto,\cpto>_X$ is the dual pairing between $X^*$ and~$X$.
The only exception from the convention for the norms is given
by the spaces $L^p$ constructed on 
$\Omega$, $\Gamma$, $Q$, and $\Sigma$, for $1\leq p\leq\infty$, 
whose norms are denoted by~$\norma\cpto_p$. 
Furthermore, we~put
\Bsist
  && H := \Ldue \,, \quad  
  V := \Huno 
  \aand
  W := \Hdue,
  \label{defspaziO}
  \\
  && \HG := \LdueG \,, \quad 
  \VG := \HunoG 
  \aand
  \WG := \HdueG,
  \label{defspaziG}
  \\
  && \calH := H \times \HG \,, \quad
  \calV := \graffe{(v,\vG) \in V \times \VG : \ \vG = v\suG}
  \non
  \\
  && \aand
  \calW := \bigl( W \times \WG \bigr) \cap \calV \,.
  \label{defspaziprod}
\Esist
In the sequel, we work in the framework of the Hilbert triplet
$(\calV,\calH,\calVp)$.
Thus, we have
$\<(g,\gG),(v,\vG)>_{\calV}=\iO gv+\iG\gG\vG$
for every $(g,\gG)\in\calH$ and $(v,\vG)\in\calV$.
Next, we introduce the generalized mean value, the related spaces and the operator $\calN$
we widely use throughout the paper.
The former is defined~by
\Beq
  \mean\gstar := \frac { \< \gstar,(1,1) >_{\calV} } { |\Omega| + |\Gamma| }
  \quad \hbox{for $\gstar\in\calVp$}
  \label{genmean}
\Eeq
and reduces to
\Beq
  \mean \pier{\gstar} = \frac { \iO v + \iG \vG } { |\Omega| + |\Gamma| }
  \quad \hbox{if $\gstar=(v,\vG)\in\calH$} \,.
  \label{usualmean}
\Eeq
Of course, the components of the pair $(1,1)$ in \eqref{genmean} 
are the constant functions~$1$ on $\Omega$ and~$\Gamma$, respectively.
We stress that the function
\Beq
  \calV \ni (v,\vG) \pier{{}\mapsto{}} \norma{\nabla v}_{\pier{H}}^2 + \norma{\juerg{\nablaG}\vG}_{\pier{\HG}}^2 + 
  |\mean(v,\vG)|^2
  \non
\Eeq
yields the square of a Hilbert norm on $\calV$ that is equivalent to the natural one.
In particular, we have, for every $(v,\vG)\in\calV$,
\Beq
  \normaVV{(v,\vG)} \leq \CO \bigl( \norma{\nabla v}_{\pier{H}} + \norma{\juerg{\nablaG}\vG}_{\pier{\HG}} + |\mean(v,\vG)| \bigr),
  \label{normaVVequiv}
\Eeq
where $\CO$ depends only on~$\Omega$.
Now, we set
\Beq
  \calVsz := \graffe{ \gstar\in\calVp : \ \mean\gstar = 0 } , \quad
  \calHz := \calH \cap \calVsz
  \aand
  \calVz := \calV \cap \calVsz .
  \label{defcalVz}
\Eeq
\Accorpa\Defspazi defspaziO defcalVz
Notice the difference between $\calVsz$ and the dual space~$\calVzp=(\calVz)^*$.
At this point, it is clear that the function
\Beq
  \calVz \ni (v,\vG) \pier{{}\mapsto{}} \norma{(v,\vG)}_{\calVz}
  := \bigl( \norma{\nabla v}_{\pier{H}}^2 + \norma{\juerg{\nablaG}\vG}_{\pier{\HG}}^2 \bigr)^{1/2}
  \label{normaVz}
\Eeq
is a Hilbert norm on $\calVz$ which is equivalent to the usual one.
This has the following consequence:
for every $\gstar\in\calVsz$, there exists a unique pair $\Xi\in\calVz$ such that
\Beq
  \iO \nabla\xi \cdot \nabla v + \iG \nablaG\xiG \cdot \nablaG\vG
  = \< \gstar , (v,\vG) >_{\calV}
  \quad \hbox{for every $(v,\vG)\in\calV$}.
  \qquad
  \label{perdefN}
\Eeq
Indeed, the \rhs\ of \eqref{perdefN}, restricted to the pairs $(v,\vG)\in\calVz$,
\juerg{defines} a continuous linear functional on~$\calVz$ with respect to its natural norm 
($\calVz$ is a subspace of~$\calV\subset V\times\VG$), \juerg{and}
thus also with respect to the norm~\eqref{normaVz}.
Therefore, by the Riesz representation theorem, there exists a unique pair $\Xi\in\calVz$ such that
\Beq
  \iO \nabla\xi \cdot \nabla v + \iG \nablaG\xiG \cdot \nablaG\vG
  = \< \gstar , (v,\vG) >_{\calV}
  \quad \hbox{for every $(v,\vG)\in\calVz$}.
  \qquad
  \non
\Eeq
On the other hand, the same relation holds true \juerg{for} $(v,\vG)=(1,1)$, since $\mean\gstar=0$.
As $\calV=\calVz\oplus\Span\{(1,1)\}$,
we obtain~\eqref{perdefN}.
This allows us to define $\calN:\calVsz\to\calVz$ by setting:
\Beq
  \hbox{for $\gstar\in\calVsz$, \ $\calN\gstar$ is the unique pair $\Xi\in\calVz$ satisfying \eqref{perdefN}}.
  \label{defN}
\Eeq
We notice that $\calN$ is linear, symmetric, \juerg{and bijective}.
Therefore, if we set
\Beq
  \norma\gstar_* := \norma{\calN\gstar}_{\calVz},
  \quad \hbox{for $\gstar\in\calVsz$},
  \label{normastar}
\Eeq
\juerg{then} we obtain a Hilbert norm on $\calVsz$, \juerg{which 
turns out to be} equivalent to the norm induced by the norm of~$\calVp$.
For a future use, we \juerg{collect} some properties of~$\calN$.
By just applying the definition, we have that
\begin{align}
  & \< \gstar , \calN\gstar >_{\calV}
  = \norma\gstar_*^2
  \qquad \hbox{if $\gstar \in \calVsz$},
  \label{propNa}
  \\[0.2cm]
  \separa
  & \iO \nabla w \cdot \nabla\xi
  + \iG \juerg{\nablaG}\wG \cdot \juerg{\nablaG}\xiG
  = \normaHH{(w,\wG)}^2
  \non
  \\
  & \quad \hbox{if $(w,\wG) \in \pier{{}\calVz{}}$ and $(\xi,\xiG) = \calN(w,\wG)$} \,.
  \label{propNb}
\end{align}
By accounting for the symmetry of~$\calN$, we also have
(\juerg{where,} here and in the sequel, $\calN$~is applied to $\calVsz$-valued functions as well)
\begin{align}
 &\pier{\< \dt\gstar , \calN\gstar >_{\calV}
  =  \frac 12 \, \frac d{dt} \, \norma\gstar_*^2
  \quad \hbox{if $\gstar \in \H1\calVsz$}}, 
  \label{propNdta}
\\[0.2cm]
    & \iO \nabla w \cdot \nabla\xi
  + \iG \juerg{\nablaG}\wG \cdot \juerg{\nablaG}\xiG
  = \frac 12 \, \frac d{dt} \, \normaHH{(w,\wG)}^2
  \non
  \\
  & \quad \hbox{if $(w,\wG)\in\L2\calV$, \ $\dt(w,\wG) \in \L2\calVsz $, \pier{\ $\Xi = \calN(\dt(w,\wG))$}} \,.
  \label{propNdtb}
\end{align}
\Accorpa\PropN propNa propNdtb

Now, we list our assumptions.
For the structure of our system, we \juerg{postulate:}
\begin{align}
  & \hbox{$\tauO$ and $\tauG$ \ are nonnegative real numbers}\pier{\, ;}
  \label{hptau}
  \\
  & \Beta,\, \BetaG : \erre \to [0,+\infty]
  \quad \hbox{are convex, proper and l.s.c.\ with} \quad
  \Beta(0) = \BetaG(0) = 0\pier{\,;}
  \qquad
  \label{hpBeta}
  \\[0.1cm]
  \separa
  & \Pi,\, \PiG : \erre \to \erre
  \quad \hbox{are of class $C^2$ with \Lip\ continuous first derivatives}.
  \qquad
  \label{hpPi}
\end{align}
We set, for convenience,
\Beq
  \beta := \partial\Beta \,, \quad
  \betaG := \partial\BetaG \,, \quad
  \pi := \Pi'
  \aand
  \piG := \PiG',
  \label{defbetapi}  
\Eeq
and assume that, \juerg{with some positive constants $C$ and $\eta$,}
\Beq
  D(\betaG) \subseteq D(\beta)
  \aand
  |\betaz(r)| \leq \eta |\betaGz(r)| + C
  \quad \hbox{for every $r\in D(\betaG)$}.
  \label{hpCC}
\Eeq
\Accorpa\HPstruttura hptau hpCC
In \eqref{hpCC}, the symbols $D(\beta)$ and $D(\betaG)$ 
denote the domains of $\beta$ and~$\betaG$, \pier{respectively.
More generally,} we use the notation $D(\calG)$ 
for every maximal monotone graph $\calG$ in $\erre\times\erre$,
as well as for the maximal monotone operators induced on $L^2$ spaces.
Moreover, for $r\in D(\calG)$,
$\calG^\circ(r)$ stands for the element of $\calG(r)$ having minimum modulus.

For the data, we \juerg{make the following assumptions:}
\begin{align}
  & u \in \L2{\Lx3}^3 , \quad
  \div u = 0 \quad \hbox{in $Q$}
  \aand
  u \cdot \nu = 0 \quad \hbox{on $\Sigma$}\pier{\,;}
  \label{hpu}
  \\
  & (\rhoz \,,\, \rhoGz) \in \calV \,, \quad
  \Beta(\rhoz) \in \Luno
  \aand
  \BetaG(\rhoGz) \in \LunoG \pier{\,;}
  \qquad
  \label{hprhoz}
  \\
  & \mz := \mean(\rhoz,\rhoGz)
  \in \mathop{\rm int}\nolimits D(\betaG) .
  \label{hpmz}
\end{align}
\Accorpa\HPdati hpu hpmz

Let us come to \juerg{our notion} of solution.
It is a \juerg{triple of pairs,} $\soluz$\juerg{, that} satisfies a rather low level of regularity, in principle.
Indeed, we just require that
\Bsist
  && \Mu \in \L2\calV ,
  \label{regmu}
  \\
  && \Rho \in \H1\calVp \cap \L\infty\calV ,
  \label{regrho}
  \\
  && \Zeta \in \L2\calH ,
  \label{regzeta}
  \\
  && \tauO\dt\rho \in \L2H
  \aand
  \tauG\dt\rhoG \in \L2\HG \,.
  \label{regdtrho}
\Esist
\Accorpa\Regsoluz regmu regdtrho
We have written, e.g., $\tauO\dt\rho$ in \eqref{regdtrho}
instead of $\dt(\tauO\rho)$.
We~similarly behave throughout the paper,
in particular in the forthcoming~\eqref{seconda},
in order to avoid a heavy notation.
The problem  to be  solved is stated in a weak form, 
\pier{owing} to the assumptions \eqref{hpu} on~$u$.
Namely, we require that
\Bsist
  && \< \dt\Rho , (v,\vG) >_{\calV}
  - \iO \rho u \cdot \nabla v
  + \iO \nabla\mu \cdot \nabla v
  + \iG \juerg{\nablaG}\muG \cdot \juerg{\nablaG}\vG
  = 0
  \non
  \\
  && \quad \hbox{\aet\ and for every $(v,\vG)\in\calV$},
  \label{prima}
  \\
  \separa
  && \tauO \iO \dt\rho \, v
  + \tauG \iG \dt\rhoG \, \vG
  + \iO \nabla\rho \cdot \nabla v
  + \iG \nablaG\rhoG \cdot \nablaG\vG
  \non
  \\
  && \quad {}
  + \iO \bigl( \zeta + \pi(\rho) \bigr) v
  + \iG \bigl( \zetaG + \piG(\rhoG) \bigr) \vG
  = \iO \mu v 
  + \iG \muG \vG
  \non
  \\
  && \quad \hbox{\aet\ and for every $(v,\vG)\in\calV$},
    \label{seconda}
  \\
  && \zeta \in \beta(\rho) \quad \aeQ
  \aand
  \zetaG \in \betaG(\rhoG) \quad \juerg{\mbox{a.e. on } \,\Sigma,}
  \label{terza}
  \\
  && \rho(0) = \rhoz
  \quad \aeO \pier{{}\aand \rhoG (0) = \rhoGz \quad \aeG .}
    \label{cauchy}
\Esist
\Accorpa\Pbl prima cauchy
\pier{We observe that any weak solution to problem \Pbl\ satisfies}
\Beq
  \dt \mean\Rho = 0,
  \quad \hbox{whence} \quad
  \mean \Rho(t) = \mz
  \quad \hbox{for every $t\in[0,T]$}.
  \label{conservation}
\Eeq
Indeed, it suffices to take $(v,\vG)=(|\Omega|+|\Gamma|)^{-1}(1,1)$ in \pier{\eqref{prima}}.

However, one can wonder whether the solution enjoys the better regularity
\Bsist
  && \dt\Rho = (\dt\rho,\dt\rhoG) \in \L2\calH
  \aand
  \Mu \in \L2\calW,
  \qquad
  \label{regmuW}
  \\
  && \Rho \in \L2\calW ,
  \label{regrhoW}
\Esist
and actually satisfies the boundary value problems presented in the Introduction,~i.e.,
\Bsist
  && \dt\rho + \nabla\rho \cdot u - \Delta\mu = 0
  \quad \aeQ,
  \label{primaO}
  \\
  && \dt\rhoG + \dn\mu - \DeltaG\muG = 0
  \quad \aeS,
  \label{primaG}
  \\
  && \tauO \dt\rho - \Delta\rho + \zeta + \pi(\rho) = \mu
  \quad \aeQ,
  \label{secondaO}
  \\
  && \tauG \dt\rhoG + \dn\rho - \DeltaG\rhoG + \zetaG + \piG(\rhoG) = \muG
  \quad \aeS.
  \label{secondaG}
\Esist
\Accorpa\Primabvp primaO primaG
\Accorpa\Secondabvp secondaO secondaG
This is not obvious.
For instance, it is not clear \pier{whether} the derivative 
$\dt\Rho$ can be replaced by $(\dt\rho,\dt\rhoG)$,
since the components of the test functions 
$(v,\vG)\in\calV$ used in \eqref{prima} are not independent.
In the first result we present, we answer the above questions.
However, for future use, it is convenient to prepare a more general tool.

\Bthm
\label{Strong}
Assume \HPstruttura\ for the structure, \eqref{hpu} for the velocity field and
\Bsist
  \soluz \in \L2{\calV\times\calV\times\calH}
  \quad \hbox{with} \quad
  (\tauO\dt\rho,\tauG\dt\rhoG) \in \L2\calH \,.
  \non
\Esist
Then, we have the following statements:\hfil\break
$i)$~if $\rho\in\pier{\L2W}$, $\dt\Rho\in\L2\calH$ and \eqref{prima} \juerg{is fulfilled},
then
\Bsist
  && \Mu \in \pier{\L1\calW}
  \quad \hbox{with} 
  \non
  \\
  && \norma\Mu_{\pier{\L1\calW}}
  \leq C_1 \bigl( \norma\Mu_{\L2\calV} + \norma{\dt\Rho}_{\L2\calH} 
  \non
  \\
  && \phantom{\norma\Mu_{\L2\calW} \leq C_1 \ } {} + \norma\rho_{\pier{\L2 W}} \norma u_{\L2H} \bigr),
  \label{muLdW}
\Esist
where $C_1$ \juerg{depends} only on~$\Omega$, and \Primabvp\ hold true as well;\hfil\break
$ii)$~if \eqref{seconda} \juerg{is satisfied},
then
\begin{align}
  & \Rho \in \L2\calW
  \quad \hbox{with} 
  \non
  \\
  & \norma\Rho_{\L2\calW}
  \leq C_2 \bigl( \norma\Rho_{\L2\calV} 
  \non
  \\
  & \phantom{\norma\Rho_{\L2\calW} \leq C_2 \ } {}
  + \norma{((\mu,\muG),(\zeta,\zetaG),(\tauO\dt\rho,\tauG\dt\rhoG))}_{\L2{\calH\times\calH\times\calH}} \bigr),\quad
  \label{rhoLdW}
\end{align}
where $C_2$ \juerg{depends} only on~$\Omega$, and \Secondabvp\ hold as well;\hfil\break
$iii)$~if $\gamma:\erre\to\erre$ is monotone and \Lip\ continuous,
and if \eqref{seconda} holds true with $\zetaG\in\gamma(\rhoG)$ \aeS,
then
\Bsist
  && \norma\zetaG_{\L2\HG}
  \leq C_3 \bigl( \norma\Rho_{\L2\calV} 
  \non
  \\
  && \phantom{\norma\zetaG_{\L2\HG} \leq C_2' \ } {}
  + \norma{((\mu,\muG),\zeta,(\tauO\dt\rho,\tauG\dt\rhoG))}_{\L2{\calH\times H\times\calH}} \bigr),
  \label{zetaGLdH}
\Esist
where $C_3$ \juerg{depends} only on~$\Omega$.

Assume, in addition, that $u$ belongs to $\L\infty{\Lx3}$ and that
\Beq
  \soluz \in \L\infty{\calV\times\calV\times\calH}
  \aand
  (\tauO\dt\rho,\tauG\dt\rhoG) \in \L\infty\calH \,.
  \non
\Eeq
Then, we have the following statements:\hfil\break
$iv)$~if  \pier{{}$\rho\in\L\infty W$,{}} $\dt\Rho\in\L\infty\calH$ and \eqref{prima} \juerg{is fulfilled,}
then
\Bsist
  && \Mu \in \L\infty\calW
  \quad \hbox{with} 
  \non
  \\
  && \norma\Mu_{\L\infty\calW}
  \leq C_4 \bigl( \norma\Mu_{\L\infty\calV} + \norma{\dt\Rho}_{\L\infty\calH}
  \non
  \\
  && \phantom{\norma\Mu_{\L\infty\calW} \leq C_3 \ } {} + \norma\rho_{\pier{\L\infty W}} \norma u_{\L\infty H} \bigr),
  \label{muLiW}
\Esist
where $C_4$ \juerg{depends} only on~$\Omega$;\hfil\break
$v)$~if \eqref{seconda} \juerg{is satisfied},
then
\begin{align}
  & \Rho \in \L\infty\calW
  \quad \hbox{with} 
  \non
  \\
  & \norma\Rho_{\L\infty\calW}
  \leq C_5 \bigl( \norma\Rho_{\L\infty\calV} 
  \non
  \\
  & \phantom{\norma\Rho_{\L\infty\calW} \leq C_5 \ } {}
  + \norma{((\mu,\muG),(\zeta,\zetaG),(\tauO\dt\rho,\tauG\dt\rhoG))}
  _{\L\infty{\calH\times\calH\times\calH}} \bigr), \quad
  \label{rhoLiW}
\end{align}
where $C_5$ \juerg{depends} only on~$\Omega$;\hfil\break
$vi)$~if $\gamma:\erre\to\erre$ is monotone and \Lip\ continuous,
and \juerg{if \eqref{seconda} holds true} with $\zetaG\in\gamma(\rhoG)$ \aeS,
then
\begin{align}
  & \norma\zetaG_{\L\infty\HG}
  \leq C_6 \bigl( \norma\Rho_{\L\infty\calV} 
  \non
  \\
  & \phantom{\norma\zetaG_{\L\infty\HG} \leq C_4' \ } {}
  + \norma{((\mu,\muG),\zeta,(\tauO\dt\rho,\tauG\dt\rhoG))}_{\L\infty{\calH\times H\times\calH}} 
  \bigr), \quad
  \label{zetaGLiH}
\end{align}
where $C_6$ \juerg{depends} only on~$\Omega$. \hfil\break
As a particular case of~$i)$ and~$ii)$, every solution to problem \Pbl\ satisfying \Regsoluz\
also fulfills \eqref{regrhoW} and \Secondabvp,
and, if $\tauO$ and $\tauG$ are strictly positive,
\eqref{regmuW} and \Primabvp\ hold \juerg{true} as well.
\Ethm

\Brem
\label{RemStrong}
We stress that all of the constants appearing in the estimates \accorpa{muLdW}{zetaGLiH}
depend only on~$\Omega$.
In particular, the constants $C_3$ and $C_6$ do not depend on~$\gamma$.
\Erem

Our next results regard the well-posedness and the continuous dependence of the solution on the velocity field. They are \juerg{as follows}:

\Bthm
\label{Wellposedness}
Assume \HPstruttura\ for the structure and \HPdati\ for the data.
Then, problem \Pbl\ has a at least \juerg{one} solution $\soluz$ 
satisfying the regularity properties \Regsoluz, \pier{\eqref{regrhoW}} and the inequality
\Bsist
  && \norma\Mu_{\pier{\L2\calV}}
  + \norma\Rho_{\H1\calVp \cap \L\infty\calV \pier{{} \cap \L2\calW}}
  \non
  \\
  && \quad {}
  + \norma\Zeta_{\L2\calH}
  + \tauO^{1/2} \norma{\dt\rho}_{\L2H}
  + \tauG^{1/2} \norma{\dt\rhoG}_{\L2\HG}
  \leq K_1, 
  \qquad
  \label{bddness}  
\Esist
for some constant $K_1$ that depends only on the structure of the system, 
$\Omega$, $T$, the initial data, and the norm of $u$ in~$\L2{\Lx3}^3$.
Furthermore, the components $\rho$ and $\rhoG$ of any solution are uniquely determined,
and the whole solution is unique if at least one of the operators $\beta$ and $\betaG$ is single-valued.
\Ethm

\Brem
\label{Bddnessbis}
By combining the statements of Theorems~\ref{Strong} and~\ref{Wellposedness},
it is clear that estimates also hold for the norms of $\Mu$ and $\Rho$ in~$\L2\calW$
with a constant $K_1'$ similar to~$K_1$.
\Erem

\Bthm
\label{Contdep}
Under the assumptions \HPstruttura\ on the structure and \HPdati\ on the data,
let $u_i$, $i=1,2$, be two choices of~$u$,
and let $\soluz$ be the difference of two corresponding solutions.
Then the inequality
\Bsist
  && \norma\Rho_{\L\infty\calVzp\cap\L2\calV}
  + \tauO^{1/2} \norma{\dt\rho}_{\L\infty H}
  + \tauG^{1/2} \norma{\dt\rhoG}_{\L\infty\HG}
  \non
  \\
  && \leq K_2 \norma{u_1-u_2}_{\L2{\Lx3}}
  \label{contdep}
\Esist
holds true for some constant $K_2$ that depends only on the structure of the system, 
$\Omega$, $T$, the initial data, and the norms of $u_i$\juerg{, $i=1,2$,} in~$\L2{\Lx3}^3$.
\Ethm

Under additional assumptions on the initial data and on the velocity~$u$,
we can ensure further regularity for the solution.
Namely, we have the following result:

\Bthm
\label{Regularity}
In addition to \juerg{the} assumptions \HPstruttura\ for the structure and \HPdati\ for the data,
suppose that $\tauO$ and $\tauG$ are strictly positive and that
\Bsist
  && u \in \H1{\Lx{3/2}} \cap \L\infty{\Lx3},
  \label{hpureg}
  \\
  && \rhoz \in \Hdue \,, \quad
  \rhoGz \in \HdueG \,, \quad
  \betaz(\rhoz) \in \Ldue
  \aand
  \betaGz(\rhoGz) \in \LdueG \,.
  \qquad 
  \label{hprhozreg}
\Esist
\Accorpa\HPdatireg hpureg hprhozreg
Then, problem \Pbl\ has a at least \juerg{one} solution $\soluz$ 
that also satisfies
\Bsist
  && \Mu \in \L\infty\calW \,, \quad
  \Rho \in \W{1,\infty}\calH \cap \H1\calV \cap \L\infty\calW
  \qquad
  \non
  \\
  && \aand \Zeta \in \L\infty\calH ,
  \label{regularity}  
  \\
  && \norma\Mu_{\juerg{L^\infty(0,T;\calW)}}
  + \norma\Rho_{\W{1,\infty}\calH\cap\H1\calV\cap\L\infty\calW}
  \non
  \\
  && \quad {}
  + \norma\Zeta_{\L\infty\calH}
  \leq K_3 ,
  \qquad
  \label{bddnessbis}
\Esist
with a constant $K_3$ that depends only on the structure of the system, 
$\Omega$, $T$, the initial data, and the norm of $u$ in~$\H1{\Lx{3/2}}\cap\L\infty{\Lx3}$.
In particular, the components $\Mu$ and $\Rho$ are bounded.
\Ethm

\Brem
\label{RemReg}
\pier{As $\Omega\subset \erre^3 $ and $\calW \subset  C^0(\overline \Omega)\times C^0(\overline\Gamma)$ due to the Sobolev inequalities,}
\juerg{from standard embedding results (cf., e.g., \cite[Sect.~8, Cor.~4]{Simon}) and \eqref{regularity} it follows
that even $\rho\in C^0(\overline Q)$ \pier{and} $\rho_\Gamma\in C^0(\overline\Sigma)$. 
Moreover, a part of the result of Theorem 2.6 still holds true}
without assuming that $\tauO$ and $\tauG$ are strictly positive,
provided that the initial data satisfy the additional condition
\Bsist
  && \bigl( {-\Delta}\rhoz + (\betaeps+\pi)(\rhoz) , {-\DeltaG}\rhoGz + \dn\rhoz + (\betaGeps+\piG)(\rhoGz) \bigl)
  \non
  \\
  && \quad \hbox{belongs to a bounded subset of~$\calV$ for every $\eps\in(0,1)$}.
  \label{hprhozregV}
\Esist
With respect to the previous statement,
we miss the conditions $\dt\Rho\in\L\infty\calH$
and $\Mu\in\L\infty\calW$
(see the forthcoming Remark~\ref{PerRemReg} for \juerg{details}).
If the double-well potentials in the bulk and on the boundary are the same potential of logarithmic type
as in the next \accorpa{stessodom}{comelogpot}, \juerg{then}
it is easy to find sufficient conditions on $\rhoz$ for \eqref{hprhozregV} to hold.
Indeed, one can assume that $\norma\rhoz_\infty<1$ and $(\Delta\rhoz,\DeltaG\rhoGz-\dn\rhoz)\in\calV$.
\Erem

Our last result requires potentials of logarithmic type (see~\eqref{logpot}) with the same domain.
Namely, we require that
\Bsist
  && \beta,\,\betaG : (-1,1) \to \erre
  \quad \hbox{are $C^2$ functions with}
  \label{stessodom}
  \\
  && \lim_{r\seto-1} \beta(r) = \lim_{r\seto-1} \betaG(r) = -\infty
  \aand
  \lim_{r\neto1} \beta(r) = \lim_{r\neto1} \betaG(r) = +\infty \,.
  \label{comelogpot}
\Esist

\Bthm
\label{Separation}
In addition to \juerg{the} assumptions \HPstruttura\ on the structure,
assume that $\tauO$ and $\tauG$ are strictly positive and that
$\beta$ and $\betaG$ satisfy \accorpa{stessodom}{comelogpot}.
Moreover, assume that $u$ and $\rho_0$ satisfy \eqref{hpu}, \eqref{hpureg} and 
\Beq
\rho_0\in W, \qquad \rho_{0|\Gamma}\in W_\Gamma, \qquad  \inf\rhoz > -1 
  \aand
  \sup\rhoz < 1 \,.
  \label{hpsepar}
\Eeq
Then the unique solution $\soluz$ satisfies 
\Beq
\juerg{
  \rhomin \leq \rho(x,t) \leq \rhomax
  \quad \hbox{for all $(x,t)\in \overline Q$,}
  }
  \label{separation}
\Eeq
for some constants $\rhomin,\rhomax\in(-1,1)$ that depend only on the structure of the system, 
$\Omega$, $T$, the initial data, and the norm of $u$ in~$\H1{\Lx{3/2}}\cap\L\infty{\Lx3}$.
\Ethm

\Bthm
\label{ContDepbis}
In addition to \HPstruttura,
assume that $\tauO$ and $\tauG$ are strictly positive, that
$\beta$ and $\betaG$ satisfy \accorpa{stessodom}{comelogpot},
and that $\beta$, $\pi$, $\betaG$ and $\piG$ are of class~$C^2$.
Moreover, assume that $\rhoz$ satisfies~\eqref{hpsepar}.
Finally, let $u_i\in\H1{\Lx3}$, $i=1,2$, be two choices of~$u$
satisfying \eqref{hpu},
and let $\soluz$ be the difference of the corresponding solutions.
Then the inequality
\Bsist
  && \norma\Mu_{\L\infty\calW}
  + \norma\Rho_{\W{1,\infty}\calH\pier{\cap\H1\calV}\cap\L\infty\calW}
  \non
  \\
  && \leq K_4 \norma{u_1-u_2}_{\H1{\Lx3}}
  \label{contdepbis}
\Esist
holds true for some constant $K_4$ that depends only on the structure of the system, 
$\Omega$, $T$, the initial data, and the norms of $u_i$, \juerg{$i=1,2$,} in~$\H1{\Lx3}$.
\Ethm

Throughout the paper, we will repeatedly use Young's inequality
\Beq
  a\,b \leq \delta\,a^2 + \frac 1{4\delta} \, b^2
  \quad \hbox{for all $a,b\in\erre$ and $\delta>0$},
  \label{young}
\Eeq
H\"older's inequality, and the Sobolev inequality 
related to the continuous embedding $V\subset L^p(\Omega)$ with $p\in[1,6]$
(since $\Omega$ is three-dimensional, bounded and smooth).
Besides, this embedding is compact for $p<6$,
and the same holds for the analogous spaces on the boundary.
It follows that the embeddings $\calV\subset\calH$ and $\calH\subset\calVp$ are compact as well.
In particular, we have the compactness inequality
\Bsist
  && \normaHH{(v,\vG)}
  \leq \delta \bigl( \normaH{\nabla v} + \normaHG{\nablaG\vG} \bigr)
  + C_\delta \, \norma{(v,\vG)}_{\calVp}
  \non
  \\
  && \quad \hbox{for every ${(v,\vG)}\in\calV$ and $\delta>0$},
  \label{compact}
\Esist
where $C_\delta $ depends only on~$\Omega$ and~$\delta$.
Finally, we set, for brevity,
\Beq
  Q_t := \Omega \times (0,t)
  \aand
  \Sigma_t := \Gamma \times (0,t)
  \quad \hbox{for $0<t\leq T$,}
  \label{defQt}
\Eeq
and simply write $Q$ and $\Sigma$ if $t=T$.

We conclude this section by stating a general rule concerning the constants 
that appear in the estimates to be performed in the sequel.
The small-case symbol $c$ stands for a generic constant
whose values might change from line to line \pier{(}and even within the same line\pier{)}
and depend only on~$\Omega$, on the shape of the nonlinearities,
and on the constants and the norms of the functions involved in the assumptions of our statements.
In particular, the values of $c$ do not depend on $\eps$ if this parameter is considered.
A~small-case symbol with a subscript like $c_\delta$ (in~particular, with $\delta=\eps$)
indicates that the constant might depend on the parameter~$\delta$, in addition.
On the contrary, we mark precise constants that we can refer~to
by using different symbols, like in~\eqref{hpCC} and in~\eqref{bddness}.


\section{Strong solutions}
\label{STRONG}
\setcounter{equation}{0}

This section is devoted to the proof of Theorem~\ref{Strong}.
Our argument relies on a result on an elliptic problem.
Thus, we prove the following lemma:

\Blem
\label{Elliptic}
Let $\gamma:\erre\to\erre$ be monotone and \Lip\ continuous,
and assume that $(w,\wG)\in\calV$ and $(g,\gG)\in\calH$ \juerg{satisfy}
\Beq
  \iO \nabla w \cdot \nabla v
  + \iG \juerg{\nablaG}\wG \cdot \juerg{\nablaG}\vG
  + \iG \gamma(\wG) \vG
  = \iO g v
  + \iG \gG \vG
  \quad \hbox{for every $(v,\vG)\in\calV$}.
  \label{elliptic}
\Eeq
Then we have that
\Beq
  (w,\wG) \in \calW
  \aand
  \normaWW{(w,\wG)}
  + \normaHG{\gamma(\wG)}
  \leq \CO \bigl( \normaVV{(w,\wG)} + \normaHH{(g,\gG)} \bigr),
  \label{ellreg}
\Eeq
where $\CO$ depends only on~$\Omega$.
Moreover, $(w,\wG)$ solves the boundary value problem
\Beq
  - \Delta w = g
  \quad \aeO,
  \aand
  \dn w - \DeltaG\wG + \gamma(\wG) = \gG
  \quad \aeG \,.
  \label{bvpbl}
\Eeq
\Elem

\proof
We use well-known estimates \juerg{from} the theory of traces and elliptic equations.
The values of $c$ \juerg{will} depend only on~$\Omega$.
We set, for brevity, $M:=\normaVV{(w,\wG)}+\normaHH{(g,\gG)}$.
By taking any $v\in\Hunoz$ and testing \eqref{elliptic} by~$(v,0)$,
we obtain the first \juerg{identity in} \eqref{bvpbl} in the sense of \juerg{distributions}.
In particular, we~have $\Delta w= \pier{{}- g}\in H$.
By combining \juerg{this} with $w\suG=\wG\in\VG$, we deduce that
\Beq
  w \in \Hx{3/2}
  \aand
  \norma w_{\Hx{3/2}}
  \leq \pier{c}\, \bigl( \normaH{\Delta w} + \normaVG\wG \bigr) 
  \leq c \, M .
  \non
\Eeq
It follows that
\Beq
  \dn w \in \HG
  \aand
  \normaHG{\dn w}
  \leq c\,  \bigl( \norma w_{\Hx{3/2}} + \normaH{\Delta w} \bigr)
  \leq c \, M, 
  \non
\Eeq
as well as the \juerg{validity of the} formula
\Beq
  \iO \nabla w \cdot \nabla v 
  = - \iO \Delta w \, v
  + \iG \dn w \, v\suG 
  \quad \hbox{for every $v\in V$}.
  \non
\Eeq
By replacing $-\Delta w$ by~$g$, comparing with~\eqref{elliptic},
and noticing that for every $\vG\in\VG$ there exists \juerg{some} $v\in V$ such that $(v,\vG)\in\calV$, 
we deduce that
\Beq
  \iG \nablaG\wG \cdot \nablaG\vG
  + \iG \gamma(\wG) \vG
  = \iG (\gG - \dn w) \vG 
  \quad \hbox{for every $\vG\in\VG$} \,.
  \label{perlemma}
\Eeq
In particular, by choosing $\vG=\gamma(\wG)$, we obtain \juerg{that}
\Beq
  \iG \gamma'(\wG) |\nablaG\wG|^2
  + \iG |\gamma(\wG)|^2
  = \iG (\gG - \dn w) \gamma(\wG),
  \non
\Eeq
whence immediately
\Beq
  \normaHG{\gamma(\wG)}
  \leq \normaHG{\gG - \dn w}  
  \leq c \, M ,
  \non
\Eeq
which is a part of~\eqref{ellreg}.
Then, we can rewrite \eqref{perlemma} in the form
\Beq
  \iG \nablaG\wG \cdot \nablaG\vG
  = \iG (\gG - \dn w - \gamma(\wG)) \vG 
  \quad \hbox{for every $\vG\in\VG\,$}.
  \non
\Eeq
This implies the second \juerg{identity in} \eqref{bvpbl} (at least in a generalized sense), 
as well~as
\Beq
  \DeltaG\wG \in \HG
  \aand
  \normaHG{\DeltaG\wG}
  \leq \normaHG{\gG - \dn w - \gamma(\wG)}
  \leq c \, M \,.
  \non
\Eeq
Therefore, we also have \juerg{that}
\Beq
  \wG \in \WG
  \aand
  \normaWG\wG
  \leq c \bigl( \normaVG\wG + \normaHG{\DeltaG\wG} \bigr) 
  \leq c \, M \,.
  \non
\Eeq
We conclude that
\Beq
  w \in W
  \aand
  \normaW w
  \leq c \bigl( \normaH{\Delta w} + \normaWG\wG \bigr)
  \leq c \, M \,.
  \non
\Eeq
Therefore, both the regularity and the estimate of \eqref{ellreg} are completely proved,
and the equations \eqref{bvpbl} hold almost everywhere.\QED

\step
Proof of Theorem~\ref{Strong}

In order to prove $i)$ and~$iv)$,
we account for \eqref{hpu}, which implies that
$-\iO\rho\,u\cdot\nabla v=\iO\nabla\rho\cdot u\,v$
\aet\ for every $v\in V$,
and rewrite \eqref{prima} \aet\ with this substitution.
Then, \aat, we apply Lemma~\ref{Elliptic} with
\Beq
  \gamma = 0 , \quad
  (w,\wG) = \Mu(t), \quad
  g = - (\dt\rho + \nabla\rho \cdot u) (t)
  \aand
  \gG = - \dt\rhoG(t),
  \non
\Eeq
by observing that $\norma{\nabla\rho(t)\cdot u(t)}_2 \pier{{}\leq \norma{\nabla\rho(t)}_6 \norma{ u(t)}_3 \leq c\normaW{\rho(t)}}\norma{u(t)}_3\,$,
where $c$ depends only on~$\Omega$.
Then, we take the norms of both sides of \eqref{ellreg}
in \pier{$L^1(0,T)$} or in $L^\infty(0,T)$ to deduce \eqref{muLdW} and~\eqref{muLiW}, respectively,
and notice that \eqref{bvpbl} coincides with \Primabvp.
To prove $ii)$ and~$v)$,
we apply Lemma~\ref{Elliptic} \aat\ with
\Bsist
  && \gamma = 0 , \quad
  (w,\wG) = \Rho(t), \quad
  g = \bigl( \mu - \tauO \dt\rho - \zeta - \pi(\rho) \bigr)(t)
  \non
  \\
  && \aand
  \gG = \bigl( \muG - \tauG \dt\rhoG - \zetaG - \piG(\rhoG) \bigr)(t),
  \non
\Esist
and argue as before.
Finally, to prove $iii)$ and~$vi)$,
we apply Lemma~\ref{Elliptic} \aat\ with $\gamma$ as in the statement,
$(w,\wG)$ and $g$ as in the previous step, and
\Beq
  \gG = \bigl( \muG - \tauG \dt\rhoG  - \piG(\rhoG) \bigr)(t).
  \non
\Eeq
Then, we write the estimate for~$\zetaG$ of \eqref{elliptic}
and take the norms of both sides in $L^2(0,T)$ or in $L^\infty(0,T)$.
\QED

\section{Continuous dependence and uniqueness}
\label{CONTDEP}
\setcounter{equation}{0}

In this section, we give the proof of Theorem~\ref{Contdep} concerning continuous dependence on the velocity field~$u$
and derive the uniqueness part of Theorem~\ref{Wellposedness}.

\step
Proof of Theorem~\ref{Contdep}

We take two choices $u_i$, $i=1,2$, of $u$ 
and consider two corresponding solutions
\juerg{$((\mu_i,\mu_{i\Gamma}),(\rho_i,\rho_{i\Gamma}),(\zeta_i,\zeta_{i\Gamma}))$}.
We set $\rho:=\rho_1-\rho_2$ and similarly define the other differences,
according to the notation of the statement.
We observe that $\mean\Rho=0$ by the conservation property \eqref{conservation}, 
applied to~$(\rho_i,\rho_{i\Gamma})$ for $i=1,2$,
whence $\Xi(s):=\calN(\Rho(s))$ is well defined for every $s\in[0,T]$.
Thus, we write equation \eqref{prima} at the time $s$ for both solutions,
test the difference by $\Xi(s)$ and integrate with respect to $s$ over~$(0,t)$,
\juerg{where} $t\in(0,T)$.
\pier{Owing} to~\eqref{propNdta}, we obtain \juerg{the identity}
\Beq
  \frac 12 \norma{\Rho(t)}_*^2
  + \intQt \nabla\mu \cdot \nabla\xi
  + \intSt \nablaG\muG \cdot \nablaG\xiG
  = \intQt (\rho_1 u_1 - \rho_2 u_2) \cdot \nabla\xi \,.
  \label{diffprima}
\Eeq
At the same time, we write equation \eqref{seconda} at the time $s$ for both solutions,
test the difference by $\Rho(s)$, integrate over~$(0,t)$,
and add the same term $\iot\normaHH{\Rho(s)}^2\,ds$ to both sides, for convenience.
We~obtain that
\Bsist
  && \frac \tauO 2 \iO |\rho(t)|^2
  + \frac \tauG 2 \iG |\rhoG(t)|^2
  + \iot \normaV{\rho(s)}^2 \, ds
  + \iot \normaVG{\rhoG(s)}^2 \, ds
  + \intQt \zeta \rho
  + \intSt \zetaG \rhoG 
  \non
  \\
  && = \intQt \bigl\{ \rho^2 - \bigl( \pi(\rho_1) - \pi(\rho_2) \bigr) \rho \bigr\}
  + \intSt \bigl\{ \rhoG^2 - \bigl( \piG(\rho_{1\Gamma}) - \piG(\rho_{2\Gamma}) \bigr) \rhoG \bigr\} 
  \non
  \\
  && \quad {}
  + \intQt \mu \rho
  + \intSt \muG \rhoG \,.
  \label{diffseconda}
\Esist
At this point, we add these equalities to each other.
By the definition of~$\calN$, the last two integrals of \eqref{diffseconda}
and the ones on the \lhs\ of \eqref{diffprima} cancel \juerg{each other}.
Moreover, the terms involving $\zeta$ and $\zetaG$ are nonnegative by monotonicity.
Thus, by owing to the \Lip\ continuity of $\pi$ and~$\piG$, 
we deduce~that
\Bsist
  && \frac 12 \norma{\Rho(t)}_*^2
  + \frac \tauO 2 \iO |\rho(t)|^2
  + \frac \tauG 2 \iG |\rhoG(t)|^2
  + \iot \normaVV{\Rho(s)}^2 \, ds
  \non
  \\
  && \leq \intQt |\rho u_1 + \rho_2 u| \, |\nabla\xi|
  + c \iot \normaHH{\Rho(s)}^2 \, ds 
  =: I_1 + I_2,
  \non
\Esist
and we now treat the contributions \pier{$I_1$ and~$I_2$ on the \rhs} separately. 
We account for the \Holder, Sobolev and Young inequalities, and use the definitions~\eqref{normaVz} and~\eqref{normastar}.
We~have that
\Bsist
  && I_1
  \leq \iot \bigl( \norma{\rho(s)}_6 \, \norma{u_1(s)}_3 + \norma{\rho_2(s)}_6 \, \norma{u(s)}_3 \bigr) \norma{\nabla\xi(s)}_2 \, ds
  \non
  \\
  \separa
  && \leq \frac 14 \iot \normaVV{\Rho(s)}^2 \, ds
  + c \iot \norma{u_1(s)}_3^2 \, \norma{\Rho(s)}_*^2 \, ds
  \non
  \\
  && \quad {}
  + c \, \norma{\rho_2}_{\L\infty V}^2 \iot \norma{u(s)}_3^2 \, ds 
  + \iot \norma{\Rho(s)}_*^2 \, ds \,.
  \non
\Esist
We deal with $I_2$ \juerg{as follows, invoking the compactness inequality~\eqref{compact}:}
\Beq
  I_2 \leq \frac 14 \iot \normaVV{\Rho(s)}^2 \, ds
  + c \iot \norma{\Rho(s)}_*^2 \, ds \,.
  \non
\Eeq
At this point, we collect all of these inequalities,
observe that the function $s\mapsto\norma{u_1(s)}_3^2$ belongs to~$L^1(0,T)$ by \eqref{hpu},
and apply the Gronwall lemma.
We immediately deduce~\eqref{contdep}
with a constant whose dependence agrees with \juerg{that asserted in} the statement of Theorem~\ref{Contdep}.
\juerg{With this, the} proof is complete.

\step
Partial uniqueness and uniqueness

\juerg{Next, we} derive the uniqueness part of Theorem~\ref{Wellposedness}.
Uniqueness for $\Rho$ clearly follows by taking $u_1=u_2$ in~\eqref{contdep}.
Assume now that $\beta$ is single-valued.
This implies that $\zeta=\beta(\rho)$ is uniquely determined as well.
Next, by Theorem~\ref{Strong}, \Secondabvp\ hold true.
From \eqref{secondaO}, we deduce uniqueness for the component $\mu$ of the solution.
This also implies uniqueness for $\muG=\mu\suS$,
and \eqref{secondaG} yields uniqueness for~$\zetaG$.
Assume now that $\betaG$ is single-valued.
In this case, we first derive uniqueness for $\zetaG=\betaG(\rhoG)$, then for $\muG$ by owing to~\eqref{secondaG}.
On the other hand, the first equation \eqref{prima} with $\Rho$ completely known
implies that the difference of the components $\Mu$ of two solutions is space independent,
whence it has the form
$t\mapsto\phi(t)(1,1)$ for some $\phi\in L^2(0,T)$, 
since the second component is the trace of the first one. 
But $\phi$ must vanish since $\muG$ is unique.
This \juerg{implies} that $\mu$ is unique as well.
Finally, \eqref{secondaO} yields uniqueness for~$\zeta$.\QED


\section{Approximation}
\label{APPROXIMATION}
\setcounter{equation}{0}

In this section, we construct and solve an approximating problem
depending on the small parameter $\eps\in(0,1)$,
which is understood to be fixed throughout the whole section.
This problem is simply obtained by modifying \Pbl\ as follows:
instead of $\tauO$ and~$\tauG$, we take the strictly positive constants
\Beq
  \tauOeps := \max \{\tauO,\eps\}
  \aand
  \tauGeps := \max \{\tauG,\eps\},
  \label{deftaueps}
\Eeq
and replace the functionals $\Beta$ and $\BetaG$ and the operators $\beta$ and $\betaG$
by \pier{the following Moreau and Yosida regularizations $\Betaeps$, $\BetaGeps$, $\betaeps$, $\betaGeps$ (see, e.g., \cite[pp.~28 and~39]{Brezis}):
\begin{align*}
	&\widehat{\beta }_{\varepsilon }(r)
	:=\inf_{s \in \mathbb{R}}\left\{ \frac{1}{2\varepsilon } |r-s|^2
	+\widehat{\beta }(s) \right\} 
	= \int_{0}^{r} \beta _\varepsilon (s)ds,
	\\
	&\widehat{\beta }_{\Gamma, \varepsilon }(r)
	:=\inf_{s \in \mathbb{R}}\left\{ \frac{1}{2\varepsilon \eta } |r-s|^2+\widehat{\beta }_\Gamma (s) \right\} 
	= \int_{0}^{r} \beta _{\Gamma,\varepsilon} (s)ds,
	\\
	&	\beta _\varepsilon (r)
	:=\frac{1}{\varepsilon }\bigl( r-(I+\varepsilon \beta )^{-1} (r)\bigr),
	\\
	&\beta _{\Gamma, \varepsilon} (r)
	:=\frac{1}{\varepsilon \eta}\bigl( r- (I+\varepsilon \eta \beta _\Gamma )^{-1} (r) \bigr)
\end{align*}
for all $ r\in \erre$, where $\eta>0$ is {\pier the} same constant as in the assumption \eqref{hpCC}. We point out that \eqref{hpBeta} and \eqref{defbetapi} hold also for the approximations. Moreover, we have that
\Bsist
  && 0 \leq \Betaeps(r) \leq \Beta(r),  \quad 
	0 \le \widehat{\beta }_{\Gamma, \varepsilon} (r) \le \widehat{\beta }_{\Gamma }(r)
  \quad \hbox{for every $r\in\erre$},
  \label{propBetaeps}
  \\
  && |\betaeps(r)| \leq |\betaz(r)| , \quad 
	\bigl |\beta _{\Gamma ,\varepsilon }(r) \bigr | \le \bigl |\beta _{\Gamma }^\circ (r) \bigr |  \quad \hbox{for every $r\in D(\beta)$}.
  \label{propbetaeps}
\Esist
Furthermore, \eqref{hpCC}~also holds true for $\betaeps$ and $\betaGeps$
with the same constants (see~\cite[Lemma~4.4]{CaCo}).
We thus write
\Beq
  |\betaeps(r)| \leq \eta |\betaGeps(r)| + C
  \quad \hbox{for every $r\in\erre$}.
  \label{propCCeps}
\Eeq
Since} $\betaeps$ and $\betaGeps$ have the same sign,
we see that \eqref{propCCeps} and the Young inequality \pier{yield}
\Beq
  \betaGeps(r) \betaeps(r)
  \geq \frac 1 {2\eta} \, |\betaeps(r)|^2 - C_\eta
  \quad \hbox{for every $r\in\erre$},
  \label{prodbetaeps}
\Eeq
with a similar constant~$C_\eta$.
We also notice that the inclusion $D(\betaG)\subseteq D(\beta)$ (see~\eqref{hpCC}) and \eqref{hpmz} imply \juerg{that}
\Beq
  \betaeps(r) (r-\mz)
  \geq \delta_0 |\betaeps(r)| - C_0
  \aand
  \betaGeps(r) (r-\mz)
  \geq \delta_0 |\betaGeps(r)| - C_0
  \label{trickMZ}
\Eeq
for every $r\in\erre$ and every $\eps\in(0,1)$,
where $\delta_0$ and $C_0$ are some positive constants that depend only on $\beta$, $\betaG$
and on the position of $\mz$ in the interior of~$D(\betaG)$ and of~$D(\beta)$
(see, e.g. \cite[p.~908]{GiMiSchi}).

The \juerg{sought solution is} a quadruple $(\mueps,\muGeps,\rhoeps,\rhoGeps)$
\juerg{having the} regularity properties
\begin{align}
  &\Mueps \in \pier{{}\L2\calV \cap \L1\calW},
  \label{regeps1} \\
  &\Rhoeps \in \H1\calH \cap \L\infty\calV \cap \L2\calW
  \label{regeps2}
\end{align}
and such that the \juerg{6-tuple} $(\mueps,\muGeps,\rhoeps,\rhoGeps,\zetaeps,\zetaGeps)$
obtained by setting
\Beq
  \zetaeps := \betaeps(\rhoeps) 
  \aand
  \zetaGeps := \betaGeps(\rhoGeps)
  \label{defzetaeps}
\Eeq
solves the following problem:
\Bsist
  && \iO \dt\rhoeps \, v + \iG \dt\rhoGeps \, \vG
  - \iO \rhoeps u \cdot \nabla v
  + \iO \nabla\mueps \cdot \nabla v
  + \iG \nabla\muGeps \cdot \nabla\vG
  = 0
  \qquad
  \non
  \\
  && \quad \hbox{\aet\ and for every $(v,\vG)\in\calV$},
  \label{primaeps}
  \\
  \separa
  && \tauOeps \iO \dt\rhoeps \, v
  + \tauGeps \iG \dt\rhoGeps \, \vG
  + \iO \nabla\rhoeps \cdot \nabla v
  + \iG \nablaG\rhoGeps \cdot \nablaG\vG
  \non
  \\
  && \quad {}
  + \iO \bigl( \zetaeps + \pi(\rhoeps) \bigr) v
  + \iG \bigl( \zetaGeps + \piG(\rhoGeps) \bigr) \vG
  = \iO \mueps v 
  + \iG \muGeps \vG
  \non
  \\
  && \quad \hbox{\aet\ and for every $(v,\vG)\in\calV$},
  \label{secondaeps}
  \\
  && \rhoeps(0) = \rhoz
  \quad \aeO \pier{{}\aand \rhoGeps (0) = \rhoGz \quad \aeG .}
  \label{cauchyeps}
\Esist
\Accorpa\Pbleps defzetaeps cauchyeps
We have written the sum of two integrals instead of a duality in \eqref{primaeps},
\juerg{in accordance with} the requirement \eqref{regeps2} on~$\Rho$.

The aim of this section is \juerg{to solve} the approximating problem \Pbleps.
\juerg{In} this respect, we have the following result.

\Bthm
\label{Wellposednesseps}
Assume \HPstruttura\ and \eqref{deftaueps} for the structure and \HPdati\ for the data.
Then the problem \Pbleps\ has a unique solution $(\mueps,\muGeps,\rhoeps,\rhoGeps)$
\juerg{with} the regularity properties~\pier{\eqref{regeps1}--\eqref{regeps2}}.
\Ethm

\medskip

The rest of the section is devoted to the proof of Theorem~\ref{Wellposednesseps}.
Since the approximating problem~\Pbleps\ is a particular case of problem \Pbl\ 
and the operators $\betaeps$ and $\betaGeps$ are single-valued,
uniqueness has been already established in the previous section.
As for existence, we use a slightly modified 
Faedo--Galerkin scheme with a proper choice of the Hilbert basis.
We introduce the operator $\calA\in\calL(\calV;\calVp)$ by setting
\Beq
  \< \calA(w,\wG) , (v,\vG) >_{\calV}
  := \iO \nabla w \cdot \nabla v
  + \iG \nablaG\wG \cdot \nablaG\vG
  \quad \hbox{for $(w,\wG),(v,\vG)\in\calV$},
  \label{defA}
\Eeq
and notice that $\calA$ is nonnegative and weakly coercive.
Indeed, we have that
\Beq
  \< \calA(v,\vG) , (v,\vG) >_{\calV}
  + \normaHH{(v,\vG)}^2 = \normaVV{(v,\vG)}^2 
  \quad \hbox{for every $(v,\vG)\in\calV$}.
  \label{coerc}
\Eeq
Moreover, as the embedding $\calV\subset\calH$ is compact,
the resolvent of $\calA$ is compact as well, and
the spectrum of $\calA$ reduces to a discrete set of eigenvalues,
the eigenvalue problem being
\Beq
  (e,e_\Gamma) \in \calV \setminus \{(0,0)\}
  \aand
  \calA(e,e_\Gamma) = \lambda(e,e_\Gamma) \,.
  \label{eigenpbl}
\Eeq
More precisely, we can rearrange the eigenvalues and choose the eigenvectors in order~that
\Bsist
  && 0 = \lambda_1 < \lambda_2 \leq \lambda_3 \leq \dots
  \aand
  \lim_{j\to\infty} \lambda_j = +\infty ,
  \label{eigenvalues}
  \\
  && \calA \Ej = \lambda_j \Ej
  \aand
  \iO \ei\ej + \iG \eGi\eGj = \delta_{ij}
  \quad \hbox{for $i,j=1,2,\dots$},
  \qquad
  \label{eigenvectors}
\Esist
and $\graffe{\Ej}$ generates a dense subspace of \juerg{both $\calV$ and  $\calH$}.
We notice that
\Beq
  \iO \nabla\ei \cdot \nabla\ej
  + \iG \nablaG\eGi \cdot \nablaG\eGj
  = \lambda_i \Bigl( \iO \ei\ej + \iG \eGi\eGj \Bigr)
  = \lambda_i \delta_{ij} 
  \quad \hbox{for $i,j=1,2,\dots$}.
  \non
\Eeq
We also observe that every element $(w,\wG)\in\calH$ can be written~as
\Beq
  (w,\wG) = \somma j1\infty w_j \Ej
  \quad \hbox{with} \quad
  \somma j1\infty |w_j|^2 = \normaHH{(w,\wG)}^2 < + \infty,
  \non
\Eeq
and that (on account of \eqref{coerc})
\Beq
  (w,\wG) \in \calV
  \quad \hbox{if and only if} \quad
  \somma j1\infty (1+\lambda_j) |w_j|^2 < + \infty \,.
  \non
\Eeq
Namely, the last sum yields the square of a norm on $\calV$
that is equivalent to~$\normaVV\cpto$. 
In particular, we have the following property 
(the~finite sum is the $\calH$-projection on the subspace $\calVn$ defined below):
\Beq
  \normaVV{(w^n,\wG^n)} \leq \CO \normaVV{(w,\wG)}
  \quad \hbox{if} \quad
  (w^n,\wG^n) = \somma j1n w_j \Ej ,
  \label{normaVproj}
\Eeq
where $\CO$ depends only on~$\Omega$.
At this point, we set
\Beq
  \calVn := \Span\graffe{\Ej:\ 1\leq j\leq n}
  \aand
  \calV_\infty := \bigcup_{j=1}^\infty \calVn
  = \Span\graffe{\Ej:\ j\geq1},
  \label{defVn}
\Eeq
and, for every~$n\geq 1$, we look for a quadruple $(\mun,\muGn,\rhon,\rhoGn)$ satisfying
\Bsist
  && \Mun \in \L2\calVn
  \aand
  \Rhon \in \H1\calVn,
  \vphantom\sum
  \label{regsoluzn}
  \\
  && \iO \dt\rhon \, v + \iG \dt\rhoGn \, \vG
  - \iO \rhon u \cdot \nabla v
  + \iO \nabla\mun \cdot \nabla v
  + \iG \pier{\nablaG\muGn \cdot \nablaG\vG}
  \non
  \\
  && \quad {}
  + \frac 1n \iO \mun v + \frac 1n \iG \muGn \vG 
  = 0
  \non
  \\
  && \quad \hbox{\aet\ and for every $(v,\vG)\in\calVn$},
  \label{priman}
  \\
  \separa
  && \tauOeps \iO \dt\rhon \, v
  + \tauGeps \iG \dt\rhoGn \, \vG
  + \iO \nabla\rhon \cdot \nabla v
  + \iG \nablaG\rhoGn \cdot \nablaG\vG
  \non
  \\
  && \quad {}
  + \iO \bigl( \betaeps(\rhon) + \pi(\rhon) \bigr) v
  + \iG \bigl( \betaGeps(\rhoGn) + \piG(\rhoGn) \bigr) \vG
  = \iO \mun v 
  + \iG \muGn \vG
  \non
  \\
  && \quad \hbox{\aet\ and for every $(v,\vG)\in\calVn$},
  \label{secondan}
  \\
  && \rhon(0) = \rhozn
  \quad \aeO, 
  \label{cauchyn}
\Esist
\Accorpa\Pbln regsoluzn cauchyn
where $\rhozn$ is defined by the conditions $(\rhozn,{\rhozn}\suG) \in \calVn$~and
\Beq
  \iO \rhozn v + \iG {\rhozn}\suG \vG
  = \iO \rhoz v + \iG \rhoGz \vG
  \quad \hbox{for every $(v,\vG)\in\calVn$}.
  \label{defrhozn}
\Eeq
Thus, $\rhozn$ is the first component of the orthogonal projection of $(\rhoz,\rhoGz)$ on~$\calVn$.
We have
\Beq
  \normaH\rhozn
  \leq \normaHH{(\rhozn,{\rhozn}\suG)}
  \leq \normaHH{(\rhoz,\rhoGz)} 
  \aand
  \normaVV{(\rhozn,{\rhozn}\suG)}
  \leq \CO \normaVV{(\rhoz,\rhoGz)}, 
  \label{stimarhozn}
\Eeq
the second one on account of \eqref{normaVproj}.

\step
The discrete problem

By \eqref{regsoluzn}, we have to look for $\Mun$ and $\Rhon$ given~by
\Beq
  \Mun(t) = \somma j1n \muj(t) \Ej
  \aand
  \Rhon(t) = \somma j1n \rhoj(t) \Ej
  \non
\Eeq
for some $\muj\in L^2(0,T)$ and $\rhoj\in H^1(0,T)$.
Let us introduce the $n$-vectors $\mubar:=(\muj)$ and $\rhobar:=(\rhoj)$.
Then, by rewriting 
\juerg{the}
system \accorpa{priman}{secondan} just with $(v,\vG)=\Ei$ for $i=1,\dots,n$,
we see that it takes the form
\Beq
  {\rhobar\,}'(t) - U(t) \, \rhobar(t) + D_n \, \mubar(t) = 0
  \aand
  B \, {\rhobar\,}'(t) + D \, \rhobar(t) + F(\rhobar(t)) = \mubar(t),
  \label{odesystem}
\Eeq
where 
$D_n:=\diag(\lambda_1+\frac 1n,\dots,\lambda_n+\frac 1n)$,
$D:=\diag(\lambda_1,\dots,\lambda_n)$, 
$F\!:\erre^n\to\erre^n$ is \Lip\ continuous,
and the matrices $U=(u_{ij})\in\L2{\erre^{n\times n}}$ and $B=(b_{ij})\in\erre^{n\times n}$ are given~by
\Beq
  u_{ij}(t) := \iO \ej u(t) \cdot \nabla\ei
  \quad \aat
  \aand
  b_{ij} := \tauOeps \iO \ej \ei + \tauGeps \iG \eGj \eGi,
  \non
\Eeq
for $i,j=1,\dots,n$.
By adding the second \juerg{identity in} 
\eqref{odesystem} to the first one \juerg{multiplied} by~$D_n^{-1}$,
we obtain the equivalent system
\Beq
  (D_n^{-1}+B) \, {\rhobar\,}'(t) + V(t) \, \rhobar(t) + F(\rhobar(t)) = 0
  \aand
  \mubar(t) = B \, {\rhobar\,}'(t) + D \, \rhobar(t) + F(\rhobar(t)),
  \non
\Eeq
where $V:=D-D_n^{-1}U$ belongs to $\L2{\erre^{n\times n}}$ and $D_n^{-1}+B$ is invertible,
as we verify.
To \juerg{this end}, we show that $B$ is positive definite.
Indeed, for any vector $y=(y_1,\dots,y_n)\in\erre^n$, by setting $(v,\vG):=\somma j1n y_j\Ej$,
we have that
\Bsist
  && (By) \cdot y = \somma {i,j}1n b_{ij} y_j y_i
  = \tauOeps \iO \somma i1n y_i \ei \somma j1n y_j \ej
  + \tauGeps \iO \somma i1n y_i \eGi \somma j1n y_j \eGj
  \non
  \\
  && = \tauOeps \iO |v|^2 + \tauGeps \iG |\vG|^2
  \geq \eps \normaHH{(v,\vG)}^2
  = \eps \norma y_{\erre^n}^2 \,.
  \non
\Esist
Hence, $D_n^{-1}+B$ is positive definite as well, thus invertible.
On the other hand, \eqref{cauchyn} is equivalent to an initial condition for~$\rhobar$.
Therefore, the discrete problem \Pbln\ has a unique solution.

At this point, our aim is to show that the solutions to the discete problem
converge to a solution to the approximating problem \Pbleps\ as $n$ tends to infinity, 
at least for a subsequence.
%
%
%
\juerg{To this end, we start estimating and find bounds that do not depend on~$n$.
On the contrary, they can depend on~$\eps$.}

\step
An a priori estimate

We test \eqref{priman}, written at the time $s$, by $\Mun(s)$ and integrate over $(0,t)$ with respect to~$s$ \juerg{to find that}
\Bsist
  && \intQt \dt\rhon \, \mun
  + \intSt \dt\rhoGn \, \muGn
  + \intQt |\nabla\mun|^2
  + \intSt |\nablaG\muGn|^2
  \non
  \\
  && \quad {}
  + \frac 1n \intQt |\mun|^2
  + \frac 1n \intSt |\muGn|^2
  = \intQt \rhon u \cdot \nabla\mun \,.
  \non
\Esist
\juerg{Next}, we test \eqref{secondan} by $\dt\Rhon(s)$, integrate over $(0,t)$ with respect to~$s$,
and add the same terms $\intQt\rhon\dt\rhon$ and $\intSt\rhoGn\dt\rhoGn$ to both sides for convenience.
We obtain \juerg{that}
\begin{align}
  & \tauOeps \intQt |\dt\rhon|^2
  + \tauGeps \intSt |\dt\rhoGn|^2
  + \frac 12 \, \normaVV{\Rhon(t)}^2
  + \iO \Betaeps(\rhon(t))
  + \iG \BetaGeps(\rhoGn(t))
  \non
  \\
  & = \frac 12 \, \normaVV{\Rhon(0)}^2
  + \iO \Betaeps(\rhon(0))
  + \iG \BetaGeps(\rhoGn(0))
  + \intQt \mun \dt\rhon
  + \intSt \muGn \dt\rhoGn 
  \non
  \\
  & \quad {}
  + \intQt \bigl( \rhon - \pi(\rhon) \bigr) \dt\rhon
  + \intSt \bigl( \rhoGn - \piG(\rhoGn) \bigr) \dt\rhoGn \,.
  \non
\end{align}
At this point, we add \juerg{these equalities}
and notice that four terms cancel.
Then, the remaining terms on the \lhs\ are nonnegative,
so that we can forget about four of them. 
Moreover, we use~\eqref{deftaueps}
and start estimating the \rhs\ (also accounting for~\eqref{normaVproj}, \eqref{propBetaeps} \pier{and \eqref{hprhoz}}).
We \juerg{then arrive at the estimate}
\begin{align}
  & \intQt |\nabla\mun|^2
  + \intSt |\nablaG\muGn|^2
  + \eps \intQt |\dt\rhon|^2
  + \eps \intSt |\dt\rhoGn|^2
  + \frac 12 \, \normaVV{\Rhon(t)}^2
  \non
  \\
  & \leq \intQt |\rhon| \, |u| \, |\nabla\mun|
  + c
  + \frac \eps 2 \intQt |\dt\rhon|^2
  + \frac \eps 2 \intSt |\dt\rhoGn|^2
  + c_\eps \intQt |\rhon|^2
  + c_\eps \intSt |\rhoGn|^2  \pier{{}+ c_\eps}\,.
  \non
\end{align}
On the other hand, the \Holder, Sobolev and Young inequalities yield \juerg{that}
\Bsist
  && \intQt |\rhon| \, |u| \, |\nabla\mun|
  \leq \iot\|\juerg{\rhon}(s)\|_6 \, \norma{u(s)}_3 \, \norma{\nabla\mun(s)}_2 \, ds
  \non
  \\
  && \leq \frac 12 \intQt|\nabla\mun|^2
  + c \iot \norma{u(s)}_3^2 \, \norma{\rhon(s)}_V^2 \, ds\,,
  \non
\Esist
and we notice that the function $s\mapsto\norma{u(s)}_3^2$ belongs to $L^1(0,T)$,
by~\eqref{hpu}.
Therefore, by rearranging and applying the Gronwall lemma, we \juerg{can infer}~that
\Beq
  \norma{\nabla\mun}_{\L2H}
  + \norma{\nablaG\muGn}_{\L2\HG}
  + \norma\Rhon_{\H1\calH\cap\L\infty\calV}
  \leq c_\eps \,.
  \label{stiman}
\Eeq

\step
Consequence

Just by \Lip\ continuity, we also have that
\Beq
  \norma{(\betaeps+\pi)(\rhon)}_{\L\infty H} 
  + \norma{(\betaGeps+\piG)(\rhoGn)}_{\L\infty\HG} \leq c_\eps \,.
  \non
\Eeq
On the other hand, if we test \eqref{secondan} by $(|\Omega|+|\Gamma|)^{-1}(1,1)$,
\juerg{then} we obtain, \aat,
\Bsist
  && |\mean\Mun(t)|
  \non
  \\
  && \leq c \,\bigl\{
    \normaH{\dt\rhon(t)}
    + \normaHG{\dt\rhoGn(t)}
    + \normaH{(\betaeps+\pi)(\rhon(t))}
    + \normaHG{(\betaGeps+\piG)(\rhoGn(t))}
  \bigr\} \,.
  \non
\Esist
Therefore, we \juerg{have shown} that $\mean\Mun$ is bounded in $L^2(0,T)$,
so that \eqref{stiman} and \eqref{normaVVequiv} allow us to conclude that
\Beq
  \norma\Mun_{\L2\calV} 
  \leq c_\eps \,.
  \label{stimanbis}
\Eeq

\step
Conclusion

We account for \accorpa{stiman}{stimanbis} 
and use standard weak and weak star compactness results,
as well as the Aubin-Lions lemma (see, e.g., \cite[Thm.~5.1, p.~58]{Lions}).
\juerg{It follows that}
\Bsist
  & \Mun \to \Mueps
  & \quad \hbox{weakly in $\L2\calV$},
  \label{convmun}
  \\
  & \Rhon \to \Rhoeps
  & \quad \hbox{weakly star in $\H1\calH\cap\pier{\L\infty\calV}$}
  \non
  \\
  && \quad \hbox{and strongly in $\L2\calH$},
  \label{convrhon}
\Esist
as $n$ tends to infinity, at least for a subsequence.
By \Lip\ continuity, we also deduce that
$(\betaeps+\pi)(\rhon)$ and $(\betaGeps+\piG)(\rhoGn)$
converge to 
$(\betaeps+\pi)(\rhoeps)$ and $(\betaGeps+\piG)(\rhoGeps)$
strongly in $\L2H$ and in $\L2\HG$, respectively.
Moreover, $\rhon u$ converges to $\rhoeps u$ weakly in $\L2\Ldue$, since $u\in\L2{\Lx3}$
and $\rhon$ is bounded in \pier{$\L\infty{\Lx6}$}, by the Sobolev inequality.
Finally, $\Rhon(0)$ converges to $\Rhoeps(0)$ at least weakly in~$\calH$,
so that \eqref{cauchyeps} is satisfied.

Now, we recall \eqref{defVn} for the definition of $\calV_\infty$,
and take an arbitrary $\calV_\infty$-valued step function~$(v,\vG)$.
Since the range of $(v,\vG)$ is finite-dimensional,
there exists \juerg{some}\, $m$\, such that
$(v,\vG)(t)\in\calV_m$ \aat.
It follows that $(v,\vG)(t)\in\calVn$ \aat\ and every $n\geq m$,
so that we can test \eqref{priman} and \eqref{secondan}, written at the time~$t$,
 by $(v,\vG)(t)$
and integrate over~$(0,T)$.
At this point, it is \sfw\ to deduce that
$\Mueps$, $\Rhoeps$ and the functions $\zetaeps$ and $\zetaGeps$ given by~\eqref{defzetaeps}
satisfy the integrated version of \accorpa{primaeps}{secondaeps} for every such step functions, namely, \juerg{we have that}
\Bsist
  && \intQ \dt\rhoeps \, v + \intS \dt\rhoGeps \, \vG
  - \intQ \rhoeps u \cdot \nabla v
  + \intQ \nabla\mueps \cdot \nabla v
  + \intS \nabla\muGeps \cdot \nabla\vG
  = 0\,,
  \qquad
  \non
  \\
  \separa
  && \tauOeps \intQ \dt\rhoeps \, v
  + \tauGeps \intS \dt\rhoGeps \, \vG
  + \intQ \nabla\rhoeps \cdot \nabla v
  + \intS \nablaG\rhoGeps \cdot \nablaG\vG
  \non
  \\
  && \quad {}
  + \intQ \bigl( \zetaeps + \pi(\rhoeps) \bigr) v
  + \intS \bigl( \zetaGeps + \piG(\rhoGeps) \bigr) \vG
  = \intQ \mueps v 
  + \intS \muGeps \vG \,.
  \non
\Esist
By density, the same equations hold true for every $(v,\vG)\in\L2\calV$.
This implies that \accorpa{primaeps}{secondaeps} hold \aet\ and for every $(v,\vG)\in\calV$,
as desired.
We notice that \eqref{primaeps} and \eqref{secondaeps} are
formally equal to \eqref{prima} and~\eqref{seconda}, respectively.
Moreover, by accounting for~\eqref{hpu}, we  can replace \juerg{the term}\,
$-\iO \rhoeps u\cdot\nabla v$\, by \juerg{the expression}\,
$\iO\nabla\rhoeps\cdot u\,v$\,
in \eqref{primaeps}
and notice that $\nabla\rhoeps\cdot u$ belongs to $\L2H$,
since $\rhoeps\in\pier{\L\infty{\Lx6}}$ and $u\in\L2{\Lx3}$.
This, and what we already know for the other terms,
allow us to apply $i)$ and $ii)$ of Theorem~\ref{Strong}.
We \juerg{then} deduce the full regularity \pier{\eqref{regeps1}--\eqref{regeps2}}, 
by starting from the lower regularity already established.\QED


\section{Existence}
\label{EXISTENCE}
\setcounter{equation}{0}

This section is devoted to the conclusion of the proof of Theorem~\ref{Wellposedness}.
Namely, we show that the solutions to the approximating problems converge
to a solution to problem \Pbl\ satisfying \eqref{bddness}.
We recall that the constant mean value property \eqref{conservation}
is \juerg{also} satisfied by the solutions to the \juerg{$\varepsilon$-}approximating problems.
In performing our estimates, we avoid the superscript $\eps$ in the notation of the solution, for simplicity,
writing it only at the end of each step.

\step
First a priori estimate

We test \eqref{primaeps} and \eqref{secondaeps}, written at the time~$s$,
by $\Mu(s)$ and $\dt\Rho(s)$, respectively.
Then, we integrate over~$(0,t)$ and sum up.
Moreover, we add the same terms \,$\intQt\rho\dt\rho$\, and \,$\intSt\rhoG\dt\rhoG$\, to both sides. Since some terms cancel \juerg{each other, we obtain the identity}
\Bsist
  && \intQt |\nabla\mu|^2
  + \intSt |\nablaG\muG|^2
  + \tauOeps \intQt |\dt\rho|^2
  + \tauGeps \intSt |\dt\rhoG|^2
  \non
  \\
  && \quad {}
  + \frac 12 \, \normaVV{\Rho(t)}^2
  + \iO \Betaeps(\rho(t))
  + \iG \BetaGeps(\rhoG(t))
  \non
  \\
  \separa
  && = \frac 12 \, \normaVV{(\rhoz,\rhoGz)}^2
  + \iO \Betaeps(\rhoz)
  + \iG \BetaGeps(\rhoGz)
  \non
  \\
  && \quad {}
  + \intQt \bigl( \rho - \pi(\rho) \bigr) \dt\rho
  + \intSt \bigl( \rhoG - \piG(\rhoG) \bigr) \dt\rhoG
  + \intQt \rho u \cdot \nabla\mu \,.
  \non
\Esist
Now, we observe that
\Beq
  \intQt \rho u \cdot \nabla\mu
  \leq \iot \norma{\rho(s)}_6 \, \norma{u(s)}_3 \, \norma{\nabla\mu(s)}_2 \, ds
  \leq \frac 12 \intQt |\nabla\mu|^2
  + c \iot \norma{u(s)}_3^2 \, \normaV{\rho(s)}^2 \, ds\,,
  \non
\Eeq
and that the function \juerg{\,$s\mapsto \norma{u(s)}_3^2\,$} belongs to $L^1(0,T)$, 
by~\eqref{hpu}.
Therefore, also on account of~\eqref{propBetaeps} and~\pier{\eqref{hprhoz}}, 
we easily conclude \juerg{from Gronwall's} lemma that
\Bsist
  && \norma{\nabla\mueps}_{\L2H}
  + \norma{\nablaG\muGeps}_{\L2\HG}
  + \norma\Rhoeps_{\L\infty\calV}
  \non
  \\
  && \quad {}
  + \norma{\Betaeps(\rhoeps)}_{\L\infty\Luno}
  + \norma{\BetaGeps(\rhoGeps)}_{\L\infty\LunoG}
  \non
  \\
  && \quad {}
  + (\tauOeps)^{1/2} \norma{\dt\rhoeps}_{\L2H}
  + (\tauGeps)^{1/2} \norma{\dt\rhoGeps}_{\L2\HG}
  \leq c \,.
  \label{primastima}
\Esist

\step
Consequence

By testing \eqref{primaeps} \juerg{with} an arbitrary $(v,\vG)\in\L2\calV$,
and owing to the assumptions \eqref{hpu} on~$u$,
we have that
\Bsist
  && \< \dt\Rho , (v,\vG) >_{\calV}
  \non
  \\
  && \leq \norma{\nabla\mu}_{\L2H} \norma v_{\L2V}
  + \norma{\nablaG\muG}_{\L2\HG} \norma\vG_{\L2\VG}
  \non
  \\
  && \quad {}
  + \norma\rho_{\L\infty{\Lx6}} \, \norma u_{\L2{\Lx3}} \, \norma{\nabla v}_{\L2\Ldue} \,.
  \non
\Esist
Then, the continuous embedding $V\subset\Lx6$ and \eqref{primastima} imply that
\Beq
  \norma{\dt\Rhoeps}_{\L2\calVp}
  \leq c \,.
  \label{stimadtrho}
\Eeq

\step
Second a priori estimate

We account for~\eqref{hpmz} and test \eqref{secondaeps} by the $\calVz$-valued function $(\rho-\mz,\rhoG-\mz)$ \aet\
without integrating with respect to time.
Setting $\alpha:=\mean\Mu$ \aet\ for a while, we obtain
\Bsist
  && \iO \betaeps(\rho) (\rho-\mz) 
  + \iG \betaGeps(\rhoG) (\rhoG-\mz)
  \non
  \\
  \separa
  && = - \, \tauOeps \iO \dt\rho (\rho-\mz)
  - \tauGeps \iG \dt\rhoG (\rhoG-\mz)
  - \iO |\nabla\rho|^2
  - \iG |\nablaG\rhoG|^2
  \non
  \\
  && \quad \juerg{{}-\iO \pi(\rho)(\rho-m_0)-\iG \pi_\Gamma(\rhoG)(\rhoG-m_0)}\non
  \\
  && \quad {}
  + \iO (\mu-\alpha) (\rho-\mz) 
  + \iG (\muG-\alpha) (\rhoG-\mz)
  \pier{\label{pier1}}
\Esist
\gianni{%
\pier{\aet . Observe that, in the  \rhs\ of \eqref{pier1}}
the integrals involving the gradients are bounded in $L^\infty(0,T)$\pier{, due to}~\eqref{primastima}.
Then, by using the inner product in~$\calH$, the corresponding Schwarz inequality,
and the \Lip\ continuity of $\pi$ and~$\piG$,
we deduce that
\begin{align*}
  & \iO \betaeps(\rho) (\rho-\mz) 
  + \iG \betaGeps(\rhoG) (\rhoG-\mz)
  \non
  \\
  \separa
  & \leq \bigl| \bigl( (\tauOeps\dt\rho,\tauGeps\dt\rhoG) , (\rho-\mz,\rhoG-\mz) \bigl)_{\calH} \bigr|
  + c
  \non
  \\
  & \quad {}
  + \bigl| \bigl( (\pi(\rho),\pi_\Gamma(\rhoG)) , (\rho-\mz,\rhoG-\mz) \bigl)_{\calH} \bigr|
  \non
  \\
  & \quad {}  + \bigl| \bigl( (\mu-\alpha,\muG-\alpha) , (\rho-\mz,\rhoG-\mz)_{\calH} \bigl) \bigr|
  \non
  \\
  \separa
  & \leq \bigl\{
    \norma{(\tauOeps\dt\rho,\tauGeps\dt\rhoG)}_{\calH}
    + c \, \norma{(\rho,\rhoG)}_{\calH}
    + c
    + \norma{(\mu-\alpha,\muG-\alpha)}_{\calH}
  \bigr\} \times {}
  \non
  \\
  & \quad {}
  \times \norma{(\rho-\mz,\rhoG-\mz)}_{\calH}
  + c \,.
  \non
\end{align*}
\pier{Hence, in view of \eqref{primastima} and}~\eqref{trickMZ}, we deduce that
\Beq
  \iO |\betaeps(\rho)|
  + \iG |\betaGeps(\rhoG)|
  \leq c \norma{(\mu-\alpha,\muG-\alpha)}_{\calH}
  + \psi_\eps
  \label{perseconda}
\Eeq
where $\psi_\eps$ is bounded in $L^2(0,T)$ uniformly with respect to~$\eps$}.
On the other hand, \pier{owing to the definition \eqref{normaVz}
and recalling} that $\norma\cpto_{\calVz}$ is a norm on $\calVz$ that is equivalent to the standard one,
we have that
\Beq
  \norma{(\mu-\alpha,\muG-\alpha)}_{\calH}
  \leq c \, \norma{(\mu-\alpha,\muG-\alpha)}_{\calVz} 
  = c \, \normaHH{(\nabla\mu,\nablaG\muG)} \,.
  \non
\Eeq
\gianni{Since the last term is bounded in $L^2(0,T)$ by \eqref{primastima}, \pier{the inequality}
\eqref{perseconda} implies that}
\Beq
  \norma{\betaeps(\rho)}_{\L2\Luno}
  + \norma{\betaGeps(\rhoG)}_{\L2\LunoG}
  \leq c \,.
  \non
\Eeq
At this point, we can test \eqref{secondaeps} by $(1,1)$ and find a bound for $\mean\Mu$ in $L^2(0,T)$.
\juerg{Combining it} with~\eqref{primastima}, we conclude that
\Beq
  \norma\Mueps_{\L2\calV}
  \leq c \,.
  \label{secondastima}
\Eeq

\step
Third a priori estimate

We test \eqref{secondaeps}, written at the time~$s$, with 
$(\betaeps(\rho),\betaeps(\rhoG))(s)$ \pier{and}
integrate over $(0,t)$ with respect to~$s$, \juerg{obtaining the identity} 
\Bsist
  && \tauOeps \iO \Betaeps(\rho(t))
  + \tauGeps \iG \Betaeps(\rhoG(t))
  + \intQt \betaeps'(\rho) |\nabla\rho|^2
  + \intS \betaGeps'(\rhoG) |\nablaG\rhoG|^2
  \non
  \\
  && \quad {}
  + \intQt |\betaeps(\rho)|^2 
  + \intSt \betaGeps(\rhoG) \, \betaeps(\rhoG)
  \non
  \\
  && = \tauOeps \iO \Betaeps(\rhoz)
  + \tauGeps \iG \Betaeps(\rhoGz)
  + \intQt \bigl( \mu - \pi(\rho) \bigr) \betaeps(\rho)
  + \intSt \bigl( \muG - \piG(\rhoG) \bigr) \betaeps(\rhoG)\,.
  \non
\Esist
All of the terms on the \lhs\ are nonnegative but the last one, 
for which we have, thanks to~\eqref{prodbetaeps},
\Beq
  \intSt \betaGeps(\rhoG) \, \betaeps(\rhoG)
  \geq \frac 1 {2\eta} \intSt |\betaeps(\rhoG)|^2 - c \,.
  \non
\Eeq
Since the \rhs\ can be easily handled 
by using the Young inequality, \eqref{propBetaeps}, \eqref{hprhoz},
 and the estimates \eqref{primastima} and~\eqref{secondastima},
we conclude that
\Beq
  \norma\zetaeps_{\L2H}
  + \norma{\betaeps(\rhoGeps)}_{\L2\HG}
  \leq c \,.
  \label{terzastima}
\Eeq

\step
Fourth a priori estimate

We apply the part~$iii)$ of Theorem~\ref{Strong}
to the solution to the approximating problem
with the choice $\gamma=\betaGeps$.
As the constant $C_3$ does not depend on~$\eps$, 
inequality \eqref{zetaGLdH} yields a bound for $\zetaG$
in terms of quantities that have already been estimated.
Hence, we conclude~that
\Beq
  \norma\zetaGeps_{\L2\HG} \leq c \,.
  \label{perquarta}
\Eeq
At this point, we can apply the part $ii)$ of Theorem~\ref{Strong}.
We thus have
\Beq
  \norma\Rhoeps_{\L2\calW}
  \leq c \,.
  \label{quartastima}
\Eeq

\step
Conclusion

We account for \accorpa{primastima}{quartastima} 
and use standard weak and weak star compactness results
as well as the Aubin-Lions lemma (see, e.g., \cite[Thm.~5.1, p.~58]{Lions}).
We have
\Bsist
  & \Mueps \to \Mu
  & \quad \hbox{weakly in $\L2\calV$},
  \non
  \\
  & \Rhoeps \to \Rho
  & \quad \hbox{weakly star in $\H1\calVp\cap\L\infty\calV\cap\L2\calW$}
  \non
  \\
  && \quad \hbox{and strongly in $\L2\calH$},
  \non
  \\
  & \tauOeps \dt\rhoeps \to \tauO\dt\rho
  & \quad \hbox{weakly in $\L2H$},
  \non
  \\
  & \tauGeps \dt\rhoGeps \to \tauG \dt\rhoG
  & \quad \hbox{weakly in $\L2\HG$},
  \non
  \\
  & (\zetaeps,\zetaGeps) \to (\zeta,\zetaG)
  & \quad \hbox{weakly in $\L2\calH$},
  \non
\Esist
as $\eps$ tends to zero, at least for a subsequence.
Moreover, $\rhoeps u$ converges to $\rho u$ weakly in $\L2\Ldue$, since $u\in\L2{\Lx3}$
and $\rhoeps$ converges to $\rho$ at least weakly star in $\L\infty{\Lx6}$.
At this point, it is \sfw\ to deduce that
$\soluz$ satisfies the integrated version of \accorpa{prima}{seconda}
with time-dependent test function $(v,\vG)\in\L2\calV$,
and this is equivalent to our formulation.
Furthermore, thanks to the strong convergence of $\Rhoeps$ to~$\Rho$
and to well-known results on \pier{maximal} monotone operators (see, e.g. \cite[\pier{Proposition~2.2}, p.~38]{Barbu}),
we derive~\eqref{terza}, i.e., $\zeta\in\beta(\rho)$ and $\zetaG\in\betaG(\rhoG)$.
Besides, $\Rhoeps(0)$ converges to $\Rho(0)$ at least weakly in~$\calVp$,
so that \eqref{cauchy} holds true as well.
Finally, the estimate~\eqref{bddness} follows from \pier{lower} semicontinuity.\QED


\section{Complements}
\label{COMPLEMENTS}
\setcounter{equation}{0}

This section is devoted to the proof of Theorems~\ref{Regularity}, \ref{Separation} and~\ref{ContDepbis}.
Our proofs rely on further a~priori estimates on the solutions to the 
\juerg{$\varepsilon$-}approximating problems.
However, in performing them, we proceed formally, for brevity.
Also in this \juerg{section}, we write the superscript $\eps$ in the notation for the solution
only at the end of each step.
From now on, we assume that $\tauO>0,\,\,\tauG>0$ and that \HPdatireg\ hold true.
We can also take $\eps\leq\min\graffe{\tauO,\tauG}$, so that
$\tauOeps=\tauO$ and $\tauGeps=\tauG$ (see~\eqref{deftaueps}).

\step
Fifth a priori estimate

We differentiate both \eqref{primaeps} and \eqref{secondaeps} with respect to time.
By noting that $\mean(\dt\Rho)=0$ by \eqref{conservation},
we test the obtained equations by $\Xi:=\calN(\dt\Rho)$ and $\dt\Rho$, respectively.
We \juerg{obtain the identities}
\Bsist
  && \intQt \dt^2\rho \, \xi
  + \intSt \dt^2\rhoG \, \xiG
  + \intQt \nabla\dt\mu \cdot \nabla\xi
  + \intSt \nablaG\dt\muG \cdot \nablaG\xiG
  \non
  \\
  && =  \pier{\intQt \dt\rho \, u \cdot \nabla\xi
  + \intQt \rho \, \dt u \cdot \nabla\xi\,,}
  \non
  \\
  \separa
  && \frac \tauO 2 \iO |\dt\rho(t)|^2
  + \frac \tauG 2 \iG |\dt\rhoG(t)|^2
  + \intQt |\nabla\dt\rho|^2
  + \intSt |\nablaG\dt\rhoG|^2
  \non
  \\
  && \quad {}
  + \intQt \betaeps'(\rho) |\dt\rho|^2
  + \intSt \betaGeps'(\rhoG) |\dt\rhoG|^2
  \non
  \\
  && = \frac \tauO 2 \iO |\dt\rho(0)|^2
  + \frac \tauG 2 \iG |\dt\rhoG(0)|^2
  \non
  \\
  && \quad {}
  - \intQt \pi'(\rho) |\dt\rho|^2
  - \intSt \piG'(\rhoG) |\dt\rhoG|^2
  + \intQt \dt\mu \dt\rho
  + \intSt \dt\muG \dt\rhoG \,.
  \non
\Esist
Now, we add these equalities to each other 
and treat the sum of the first two integrals by accounting for \eqref{propNdta}.
Moreover, we \juerg{can} cancel four terms in the sum \pier{due} to the definition of~$\calN$ 
(see \accorpa{perdefN}{defN}).
Finally, we recall that $\betaeps'$ and $\betaGeps'$ are nonnegative,
and integrate by parts the integrals involving~$u$ by using~\eqref{hpu}.
We then obtain that
\Bsist
  && \frac 12 \, \norma{\dt\Rho(t)}_*^2
  + \frac \tauO 2 \iO |\dt\rho(t)|^2
  + \frac \tauG 2 \iG |\dt\rhoG(t)|^2
  + \intQt |\nabla\dt\rho|^2
  + \intSt |\nablaG\dt\rhoG|^2
  \quad
  \non
  \\
  && \leq I_0
  \pier{{}- \intQt \nabla\dt\rho \cdot u \, \xi
  - \intQt \nabla\rho \cdot \dt u \, \xi{}}
  - \intQt \pi'(\rho) |\dt\rho|^2
  - \intSt \piG'(\rhoG) |\dt\rhoG|^2\,,
  \label{perquinta}
\Esist
where 
\Beq
  I_0 := \frac 12 \, \norma{\dt\Rho(0)}_*^2
  + \frac \tauO 2 \iO |\dt\rho(0)|^2
  + \frac \tauG 2 \iG |\dt\rhoG(0)|^2 \,.
  \label{defIz}
\Eeq
Now, we estimate the integrals involving~$u$
by using the \Holder\ inequality, the continuous embedding $V\subset\Lx6$,
the equivalence on $\calVz$ of the norms $\normaVV\cpto$ and~$\norma\cpto_{\calVz}$,
and the definition \eqref{normastar} of $\norma\cpto_*$.
We have 
\Bsist
  && \pier{{}-{}} \intQt \nabla\dt\rho \cdot u \, \xi
  \leq \iot \norma{\nabla\dt\rho(s)}_2 \, \norma{u(s)}_3 \, \norma{\xi(s)}_6 \, ds
  \non
  \\
  && \leq \frac 12 \intQt |\nabla\dt\rho|^2
  + c \iot \norma{u(s)}_3^2 \, \normaV{\xi(s)}^2 \, ds
  \non
  \\
  && \leq \frac 12 \intQt |\nabla\dt\rho|^2
  + c \iot \norma{u(s)}_3^2 \, \norma{\Xi(s)}_{\calVz}^2 \, ds
  \non
  \\
  && \leq \frac 12 \intQt |\nabla\dt\rho|^2
  + c \iot \norma{u(s)}_3^2 \, \norma{\dt\Rho(s)}_*^2 \, ds\,,
  \non
\Esist
as well as
\Bsist
  && \pier{{}-{}} \intQt \nabla\rho \cdot \dt u \, \xi 
  \leq \iot \norma{\nabla\rho(s)}_6 \, \norma{\dt u(s)}_{3/2} \, \norma{\xi(s)}_6 \, ds
  \non
  \\
  && \leq c \iot \normaV{\nabla\rho(s)}^2\,\juerg{ds}
  + c \iot \norma{\dt u(s)}_{3/2}^2 \, \norma{\dt\Rho(s)}_*^2 \, ds\,,
  \non
\Esist
and we notice \juerg{that} \pier{the first term on the \rhs\ is already bounded due to \eqref{quartastima}. In addition,}
the functions $s\mapsto\norma{u(s)}_3^2$ and $s\mapsto\norma{\dt u(s)}_{3/2}^2$ 
belong to $L^1(0,T)$, by \eqref{hpu} and~\eqref{hpureg}.
The last two terms on the \rhs\ of \eqref{perquinta} 
can \juerg{easily be} dealt with, by using the boundedness of \juerg{$\pi'$} and~$\piG'$
and the compactness inequality~\eqref{compact} \juerg{in the
following} way:
\Bsist  
  && - \intQt \pi'(\rho) |\dt\rho|^2
  - \intSt \piG'(\rhoG) |\dt\rhoG|^2
  \non
  \\
  && \leq \frac 12 \intQt |\nabla\dt\rho|^2
  + \frac 12 \intSt |\nablaG\dt\rhoG|^2 
  + c \iot \norma{\dt\Rho(s)}_*^2 \, ds \,.
  \non
\Esist
It remains to estimate the terms appearing in \eqref{defIz}.
To do that, we write \accorpa{primaeps}{secondaeps} at time $t=0$ \juerg{and 
account for the initial condition~\eqref{cauchyeps}. We have}
\begin{align}
  & \iO \dt\rho(0) v
  + \iG \dt\rhoG(0) \vG
  + \iO \nabla\mu(0) \cdot \nabla v
  + \iG \nablaG\muG(0) \cdot \nablaG\vG
  = \pier{ \iO \rhoz \, u(0) \cdot \nabla v\,,}
  \non
  \\
  & \tauO \iO \dt\rho(0) \, v
  + \tauG \iG \dt\rhoG(0) \, \vG
  + \iO \nabla\rhoz \cdot \nabla v
  + \iG \nablaG\rhoGz \cdot \nablaG\vG
  \non
  \\
  & \quad {}
  + \iO (\betaeps + \pi)(\rhoz)  v
  + \iG (\betaGeps + \piG)(\rhoGz) \vG
  = \iO \mu(0) v 
  + \iG \muG(0) \vG \,,
  \non
\end{align}
for every $(v,\vG)\in\calV$.
Now, we choose $(v,\vG)=\Xi:=\calN(\dt\Rho(0))$ in the first equality, 
$(v,\vG)=\dt\Rho(0)$ in the second, and add.
The terms involving $\mu(0)$ and $\muG(0)$ cancel out by the definition of~$\calN$
(see \accorpa{perdefN}{defN}).
Moreover, \juerg{invoking ~\eqref{propNa}, we obtain that}
\Bsist
  && \norma{\dt\Rho(0)}_*^2
  + \tauO \iO |\dt\rho(0)|^2
  + \tauG \iG |\dt\rhoG(0)|^2
  \non
  \\
  && = \pier{{} \iO \rhoz \, u(0) \cdot \nabla\xi {}}
  - \iO \nabla\rhoz \cdot \nabla\dt\rho(0)
  - \iG \nablaG\rhoGz \cdot \nablaG\dt\rhoG(0)
  \non
  \\
  && \quad {}
  - \iO (\betaeps + \pi)(\rhoz) \, \dt\rho(0)
  - \iG (\betaGeps + \piG)(\rhoGz) \, \dt\rhoG(0)\,, 
  \non
\Esist
and we start estimating the \rhs.
For the first term, we account for
the equivalence on $\calVz$ of the norms $\normaVV\cpto$ and~$\norma\cpto_{\calVz}$,
and the definition \eqref{normastar} of $\norma\cpto_*$ once more.
\gianni{Furthermore, we use the continuous embedding $W=\Hdue\subset\Cx0$
and the interpolation property,
where $p,p_0,p_1\in[1,+\infty]$ and $\badtheta\in(0,1)$ 
satisfy $p_0\not=p_1$ and $\frac1p=\frac{1-\badtheta}{p_0}+\frac\badtheta{p_1}$ 
(see \cite[p.~8 and Thm.~5.3.1 p.~113]{BL}),
\Beq
  (\Lx{p_0},\Lx{p_1})_{\badtheta,p}
  = (L_{p_0p_0}(\Omega),L_{p_1p_1}(\Omega))_{\badtheta,p}
  = L_{pp}(\Omega)
  = \Lx p 
  \non
\Eeq
which gives in particular $(\Lx3,\Lx{3/2})_{1/2,2}=\Ldue$ and thus the inequality
\Beq
  \norma{u(0)}_2 \leq c \, \norma u_{\H1{\Lx{3/2}}\cap\L2{\Lx3}}
  \leq c \,.
  \non
\Eeq
Hence, we can do the following computation:}
\Bsist
  && - \iO \rhoz \, u(0) \cdot \nabla\xi
  \leq \norma\rhoz_\infty \, \norma{u(0)}_2 \, \norma{\nabla\xi}_2
  \non
  \\
  && \leq c \normaW\rhoz \, \norma\Xi_{\calVz} 
  \leq c \, \norma{\dt\Rho(0)}_*
  \leq \frac 12 \, \norma{\dt\Rho(0)}_*^2 + c \,.
  \non
\Esist
We deal with the next two integrals \juerg{by} integrating by parts and using some of the assumptions~\eqref{hprhozreg}:
\Bsist
  && - \iO \nabla\rhoz \cdot \nabla\dt\rho(0)
  - \iG \nablaG\rhoGz \cdot \nablaG\dt\rhoG(0)
  \non
  \\
  && = \iO \Delta\rhoz \, \dt\rho(0)
  + \iG (\DeltaG\rhoGz - \dn\rhoz) \dt\rhoG(0)
  \leq \delta \iO |\dt\rho(0)|^2
  + \delta \iG |\dt\rhoG(0)|^2
  + c_\delta ,
  \non
\Esist
where $\delta>0$ is arbitrary.
By \juerg{invoking} \eqref{propbetaeps} for $\betaeps$ and $\betaGeps$,
\pier{and the assumptions}~\eqref{hprhozreg},
which also imply boundedness for $\rhoz$ and~$\rhoGz$, 
we \juerg{find that}
\Bsist
  && - \iO (\betaeps + \pi)(\rhoz) \, \dt\rho(0)
  - \iG (\betaGeps + \piG)(\rhoGz) \, \dt\rhoG(0) 
  \non
  \\
  && \leq \bigl( \norma{\betaz(\rhoz)}_2 + c \bigr) \norma{\dt\rho(0)}_2
  + \bigl( \norma{\betaGz(\rhoGz)}_2 +  c \bigr) \norma{\dt\rhoG(0)}_2
  \non
  \\
  && \leq \delta \norma{\dt\rho(0)}_2^2
  + \delta \norma{\dt\rhoG(0)}_2^2
  + c_\delta \,.
  \non
\Esist
Recalling all of the \juerg{above} estimates, and choosing $\delta>0$ small enough,
we see that $I_0\leq c$.
At this point, we come back to \eqref{perquinta} and apply the Gronwall lemma.
We then conclude that
\Beq
  \norma{\dt\Rhoeps}_{\L\infty\calH\cap\L2\calV} \leq c \,,
  \quad \hbox{whence} \quad
  \norma\Rhoeps_{\W{1,\infty}\calH\cap\H1\calV} \leq c \,.
  \label{quintastima}
\Eeq

\Brem
\label{PerRemReg}
In connection with Remark~\ref{RemReg},
if $\tauO$ and $\tauG$ are not supposed to be positive and \eqref{hprhozregV} holds,
one modifies the last estimates on the initial values as follows: we have
\Bsist
  && - \iO \nabla\rhoz \cdot \nabla\dt\rho(0)
  - \iG \nablaG\rhoGz \cdot \nablaG\dt\rhoG(0)
  \non
  \\
  && \quad {}
  - \iO (\betaeps + \pi)(\rhoz) \, \dt\rho(0)
  - \iG (\betaGeps + \piG)(\rhoGz) \, \dt\rhoG(0) 
  \non
  \\
  && = - \iO \bigl( {-\Delta}\rhoz + (\betaeps + \pi)(\rhoz) \bigr) \dt\rho(0)
  - \iG \bigl( \DeltaG\rhoGz - \dn\rhoz + (\betaGeps + \piG)(\rhoGz) \bigr) \dt\rhoG(0) 
  \non
  \\
  && \leq \normaVV{{-\Delta}\rhoz + (\betaeps + \pi)(\rhoz),\DeltaG\rhoGz - \dn\rhoz + (\betaGeps + \piG)(\rhoGz)}
    \, \norma{\dt\Rho(0)}_{\calVp}
  \non
  \\
  && \leq \delta \norma{\dt\Rho(0)}_{\calVp}^2
  + c_\delta \,.
  \non
\Esist
This leads to an estimate that is somewhat weaker than~\eqref{quintastima}
and yields a weaker result at the end of the procedure, as announced in the quoted remark.
\Erem

\step
Sixth a priori estimate

We set $\alpha:=\mean\Mu$ for a while
and test \juerg{\eqref{primaeps}} by the $\calVz$-valued function $\Mu-\alpha(1,1)$.
We obtain\juerg{, for a.e. $t\in (0,T)$,}
\Beq
  \iO |\nabla\mu|^2
  + \iG |\nablaG\muG|^2
  = - \iO \dt\rho (\mu - \alpha)
  - \iG \dt\rhoG (\muG - \alpha)
  \pier{{}+ \iO \rho \, u \, \nabla\mu} \,.
  \non
\Eeq
Now, we recall that the norm \eqref{normaVz} is equivalent on $\calVz$ to the natural norm.
\juerg{Thus, by} also accounting for \eqref{hpureg} 
and for \eqref{quintastima}, combined with the continuous embedding $V\subset\Lx6$,
we \juerg{may} estimate the \rhs\ \aet\ \juerg{as follows:}
\Bsist
  && - \iO \dt\rho (\mu - \alpha)
  - \iG \dt\rhoG (\muG - \alpha)
  \pier{{}+\iO \rho \, u \, \nabla\mu} 
  \non
  \\
  && \leq \pier{c\,\norma{\dt\Rho}_{\calV^*}} \norma{\Mu-\alpha(1,1)}_{\calVz}
  + \norma\rho_6 \, \norma u_3 \, \norma{\nabla\mu}_2
  \leq c \bigl( \norma{\nabla\mu}_2 + \norma{\nablaG\muG}_2 \bigr).
  \non
\Esist
At this point, the Young inequality \juerg{immediately yields that}
\Bsist
  && \norma{\nabla\mueps}_{\L\infty H} + \norma{\nablaG\muGeps}_{\L\infty\HG} \leq c \,,
  \quad \hbox{i.e.,} \quad
  \non
  \\
  && \norma{\Mueps-\mean\Mueps}_{\L\infty\calV} \leq c \,.
  \label{sestastima}
\Esist

\step
Seventh a priori estimate

We recall the estimate \eqref{perseconda} already obtained,
which holds \aet\ and also involves $\alpha:=\mean\Mu$.
\juerg{From} \eqref{quintastima} and \eqref{sestastima}, we infer that
\Beq
  \norma{\betaeps(\rho)}_{\L\infty\Luno}
  + \norma{\betaGeps(\rhoG)}_{\L\infty\LunoG}
  \leq c.
  \non
\Eeq
\juerg{We use this bound} and \eqref{quintastima} in the next estimate:
we test \eqref{secondaeps} by $(1,1)/(|\Omega|+|\Gamma|)$ and 
obtain,  \aat,
\Bsist
  && |\mean\Mu(t)|
  \leq c \, \norma{\dt\Rho}_{\pier{\L\infty{\calV^*}}}
  \non
  \\
  && \quad {}
  + c \norma{(\betaeps+\pi)(\rho)}_{\L\infty\Luno}
  + c \norma{(\betaGeps+\piG)(\rhoG)}_{\L\infty\LunoG}
  \leq c \,.
  \non
\Esist
\juerg{Combining this} with \eqref{sestastima}, we conclude that
\Beq
  \norma\Mueps_{\L\infty\calV} \leq c, 
  \quad \hbox{whence} \quad
  \norma\Mueps_{\L\infty\calH} \leq c \,.
  \label{settimastima}
\Eeq

\step
\pier{Eighth} estimate

At this point, we can test \eqref{secondaeps} by $(\betaeps(\rho),\betaeps(\rhoG))$ \aet.
By taking advantage of the above estimates and of~\eqref{prodbetaeps}, we immediately deduce that
\Beq
  \norma{\betaeps(\rhoeps)}_{\L\infty H} + \norma{\betaeps(\rhoGeps)}_{\L\infty\HG} \leq c \,.
  \label{ottavastima}
\Eeq

\step
Ninth a priori estimate

We apply the part~$vi)$ of Theorem~\ref{Strong}
to the solution to the approximating problem
with the choice $\gamma=\betaGeps$.
As the constant $C_6$ does not depend on~$\eps$, 
inequality \eqref{zetaGLiH} yields a bound for $\zetaG$
in terms of quantities that have already been estimated.
Hence, we conclude~that
\Beq
  \norma\zetaGeps_{\L\infty\HG} \leq c \,.
  \label{pernona}
\Eeq
At this point, we can apply the part $v)$ of Theorem~\ref{Strong}.
We thus have
\Beq
  \norma\Rhoeps_{\L\infty\calW}
  \leq c \,.
  \label{nonastima}
\Eeq

\step
Proof of Theorem~\ref{Regularity}

We come back to the argument used for the existence part
of proof of Theorem~\ref{Wellposedness},
recalling that the solution to the approximating problem
converges to a solution to problem \Pbl\ in a proper topology,
at least for a subsequence.
\juerg{In view of the} estimates \accorpa{quintastima}{nonastima},
the limiting solution also satisfies
the further regularity specified by \eqref{regularity},
and estimate \eqref{bddnessbis} follows from semicontinuity.\QED

\step
Proof of Theorem~\ref{Separation}

We recall that $\mu$ and $\muG$ are bounded by Theorem~\ref{Regularity}\juerg{.
Thus, accounting for \eqref{comelogpot} and \eqref{hpsepar}, 
we may} choose $\rhomin,\rhomax\in(-1,1)$ 
with $\rhomin \leq \rhoz \leq \rhomax$ 
such that
\Bsist
  && (\beta+\pi)(r) + \norma\mu_\infty \leq 0
  \aand
  (\betaG+\piG)(r) + \norma\muG_\infty \leq 0
  \quad \hbox{for every $r\in(-1,\rhomin)$},
  \non
  \\
  && (\beta+\pi)(r) - \norma\mu_\infty \geq 0
  \aand
  (\betaG+\piG)(r) - \norma\muG_\infty \geq 0
  \quad \hbox{for every $r\in(\rhomax,1)$}.
  \non
\Esist
Then, we test \eqref{seconda} by $((\rho-\rhomax)^+,(\rhoG-\rhomax)^+)$,
where $(\cpto)^+$ stands for the positive part,
and integrate with respect to time.
We obtain \juerg{the identity}
\Bsist
  && \tauO \iO |(\rho(t)-\rhomax)^+|^2
  + \tauG \iG |(\rhoG(t)-\rhomax)^+|^2
  \non
  \\
  && \quad {}
  + \intQt |\nabla(\rho-\rhomax)^+|^2
  + \intSt |\nablaG(\rhoG-\rhomax)^+|^2
  \non
  \\
  && = \intQt \bigl( \mu - (\beta+\pi)(\rho) \bigr) (\rho-\rhomax)^+
  \non
  \\
  && \quad {}
  + \intSt \bigl( \muG - (\betaG+\piG)(\rhoG) \bigr) (\rhoG-\rhomax)^+ \,.
  \non
\Esist
Since the \rhs\ is nonpositive, we conclude that $(\rho-\rhomax)^+=0$, i.e., $\rho\leq\rhomax$.
In the same way, one proves that $(\rhomin-\rho)^+=0$, i.e., $\rho\geq\rhomin$.\QED

\medskip

Now, we start the proof of Theorem~\ref{ContDepbis}.
Also in this case, we proceed formally.
Moreover, in order to simplify the notation,
we perform our estimates on the solutions to problem \Pbl, directly,
and avoid the approximating problem.
For $i=1,2$, we denote, by $\mu_i$, $\mu_{i\Gamma}$, etc., 
the components of the solutions corresponding to~$u_i$,
while $\mu$, $\muG$, etc., are the differences, e.g., $\mu=\mu_1-\mu_2$, 
according to the notation of the statement.
\juerg{For brevity, we also set $u:=u_1-u_2$,} as well~as
\Beq
  f := \Beta + \Pi \,, \quad
  \fG := \BetaG + \PiG \,,
  \quad \hbox{whence} \quad
  f' = \beta + \pi
  \aand
  \fG' = \betaG + \piG \,.
  \non
\Eeq
Moreover, since the result given by Theorem~\ref{Separation} holds for both \juerg{solutions},
we can assume that $f'$, $f''$, $\fG'$ and $\fG''$ are bounded and \Lip\ continuous,
the corresponding constants depending only on the previous assumptions on the structure,
the norms of the velocity fields $u_i$ related to~\eqref{hpureg},
and the assumptions \eqref{hprhozreg} on the initial datum.

\step
First auxiliary estimate

We write \eqref{prima} for both solutions,
take the difference and differentiate with respect to time.
Then, we test the obtained equality by $\Xi:=\calN(\dt\Rho)$ \aet\ and integrate over~$(0,t)$.
\pier{With the help of}~\eqref{propNdta} and~\eqref{normastar} we \pier{infer} that
\Beq
  \frac 12 \, \norma{\dt\Rho(t)}_*^2
  + \intQt \nabla\dt\mu \cdot \nabla\xi
  + \intSt \nablaG\dt\muG \cdot \nablaG\xiG
  = \intQt \dt \bigl( \rho_1 u_1 - \rho_2 u_2 \bigr) \cdot \nabla\xi \,.
  \non
\Eeq
At the same time, we write \eqref{seconda} for both solutions,
\pier{take the difference and differentiate it} with respect to time\pier{; then, we}
test by $\dt\Rho$ and integrate over~$(0,t)$.
Finally, we add the same integrals 
$\intQt(\rho\,\dt\rho+\nabla\rho\cdot\nabla\dt\rho)$ 
and $\intSt(\rhoG\dt\rhoG+\nablaG\rhoG\cdot\nablaG\dt\rhoG)$ 
to both sides, for convenience.
We obtain that
\Bsist
  && \frac \tauO 2 \iO |\dt\rho(t)|^2
  + \frac \tauG 2 \iG |\dt\rhoG(t)|^2
  + \intQt |\nabla\dt\rho|^2
  + \intSt |\nablaG\dt\rhoG|^2
  \non
  \\
  && \quad {}
  + \frac 12 \, \normaV{\rho(t)}^2
  + \frac 12 \, \normaVG{\rhoG(t)}^2
  \non
  \\
  \separa
  && = - \intQt \bigl( f''(\rho_1) \dt\rho_1 - f''(\rho_2) \dt\rho_2 \bigr) \dt\rho
   - \intSt \bigl( \fG''(\rho_{1\Gamma}) \dt\rho_{1\Gamma} - \fG''(\rho_{2\Gamma}) \dt\rho_{2\Gamma} \bigr) \dt\rhoG 
   \non
   \\
   && \quad {}
   + \intQt \dt\mu \, \dt\rho
   + \intSt \dt\muG \, \dt\rhoG 
   \non
   \\
   && \quad {}
   + \intQt \bigl( \rho \, \dt\rho + \nabla\rho \cdot \nabla\dt\rho \bigr)
   + \intSt \bigl( \rhoG \dt\rhoG + \nablaG\rhoG \cdot \nablaG\dt\rhoG \bigr) \,.
   \non
\Esist
At this point, we add these equalities to each other
and \juerg{employ} the definition of~$\calN$
(see \accorpa{perdefN}{defN})
in~order to cancel four terms in the sum.
Moreover, we rearrange the \rhs, 
account for \pier{\eqref{normastar}} and the equivalence of \eqref{normaVz}
to the norm in $\calV$ on the subspace~$\calVz$,
and use the boundedness and the \Lip\ continuity of both $f''$ and~$\fG''$.
We then obtain that
\begin{align}
  & \iO |\nabla\xi(t)|^2
  + \iG |\nablaG\xiG(t)|^2
  + \iO |\dt\rho(t)|^2
  + \iG |\dt\rhoG(t)|^2
  \non
  \\
  & \quad {}
  \pier{{}+ \intQt |\nabla\dt\rho|^2
  + \intSt |\nablaG\dt\rhoG|^2}
  + \normaV{\rho(t)}^2
  + \normaVG{\rhoG(t)}^2
  \non
  \\
  & \leq c \intQt |\dt\rho| \, |u_1| \, |\nabla\xi|
  + c \intQt |\dt\rho_2| \, |u| \, |\nabla\xi|
  + c \intQt |\rho| \, |\dt u_1| \, |\nabla\xi|
  + c \intQt |\rho_2| \, |\dt u| \, |\nabla\xi|
  \non
  \\
  & \quad {}
  + c \intQt |\rho| \, |\dt\rho_1| \, |\dt\rho|
  + c \intSt |\rhoG| \, |\dt\rho_{1\Gamma}| \, |\dt\rhoG|
  \non
  \\
  & \quad {}\pier{{}+ c \intQt \tonde{ |\rho_2| +1} |\dt\rho|^2
  + c \intSt \tonde{|\rho_{2\Gamma}| +1}  |\dt\rhoG|^2}
  \non
  \\
  & \quad {}
  + c \iot \normaV{\rho(s)} \, \normaV{\dt\rho(s)} \, ds
  + c \iot \normaVG{\rhoG(s)} \, \normaVG{\dt\rhoG(s)} \, ds
  \leq c \somma j1{\pier{10}} I_j\,,
  \non
\end{align}
with obvious definitions of~$I_1,\dots,\pier{I_{10}}$.
We now estimate each of \juerg{these} integrals
by using the \Holder, Sobolev and Young inequalities as follows.
We have, for every $\delta>0$,
\Bsist
  && I_1 
  \leq \iot \norma{\dt\rho(s)}_6 \, \norma{u_1(s)}_3 \, \norma{\nabla\xi(s)}_2 \, ds
  \non
  \\
  &&\quad \leq \delta \iot \normaV{\dt\rho(s)}^2 \, ds
  + c_\delta \iot \norma{u_1(s)}_3^2 \, \norma{\nabla\xi(s)}_2^2 \, ds\,, 
  \non
  \\
  \separa
  && I_2 
  \leq \iot \norma{\dt\rho_2(s)}_6 \, \norma{u(s)}_3 \, \norma{\nabla\xi(s)}_2 \, ds
  \non
  \\
  && \quad\leq \delta \iot \norma{u(s)}_3^2 \, ds
  + c_\delta \iot \normaV{\dt\rho_2(s)}^2 \, \norma{\nabla\xi(s)}_2^2 \, ds\,, 
  \non
  \\
  \separa
  && I_3 
  \leq \iot \norma{\rho(s)}_6 \, \norma{\dt u_1(s)}_3 \, \norma{\nabla\xi(s)}_2 \, ds
  \non
  \\
  && \quad\leq  \juerg{c}\iot \normaV{\rho(s)}^2 \, ds
  + \iot \norma{\dt u_1(s)}_3^2 \, \norma{\nabla\xi(s)}_2^2 \, ds\,, 
  \non
  \\
  \separa
  && I_4 
  \leq  \pier{\intQt |\nabla\xi(s)|^2
  + c {}}\, \juerg{\norma{\rho_2}_\infty^2} \intQt |\dt u|^2\,, 
  \non
  \\
  \separa
  && I_5
  \leq \iot \norma{\rho(s)}_3 \, \norma{\dt\rho_1(s)}_3 \, \norma{\dt\rho(s)}_3 \, ds
  \non
  \\
  && \quad\leq \delta \iot \normaV{\dt\rho(s)}^2 \, ds
  + c_\delta \iot \normaV{\dt\rho_1(s)}^2 \, \normaV{\rho(s)}^2 \, ds
  \non
  \\
  && \quad = \, \delta \intQt |\nabla\dt\rho|^2
  + \delta \intQt |\dt\rho|^2
  + c_\delta \iot \normaV{\dt\rho_1(s)}^2 \, \normaV{\rho(s)}^2 \, ds \,.
  \non
\Esist
Moreover, an~analogous estimate holds for~$I_6$. \pier{On the other hand, it is easy to see that
$$ I_7 \leq c\tonde{1+ \norma{\rho_2}_\infty} \intQt |\dt \rho|^2 \, , 
\quad \  I_8 \leq c \tonde{1+ \norma{\rho_{2\Gamma}}_\infty} \intSt |\dt\rhoG|^2 \, .
$$
Finally, $I_9$~and $I_{10}$} can be treated just with the Young inequality.
Now, we observe that the functions
\Beq
  s \mapsto \norma{u_1(s)}_3^2 \,, \quad 
  s \mapsto \pier{\norma{\dt\rho_i(s)}_V^2 \,, \ \, i=1,2 \, ,} \quad 
  s \mapsto \norma{\dt u_1(s)}_3^2\,, \quad
  s \mapsto \pier{\norma{\dt\rho_{1\Gamma}(s)}_{\VG}^2}\, , 
  \non
\Eeq
all \juerg{belong} to $L^1(0,T)$\pier{. Hence, we}
collect all the inequalities we have obtained,
choose $\delta$ small enough, and apply the Gronwall lemma.
We conclude~that
\Beq
  \norma{(\nabla\xi,\nablaG\xiG)}_{\L\infty\calH}
  + \norma\Rho_{\W{1,\infty}\calH\cap\H1\calV}
  \leq c \, \norma u_{\H1{\Lx3}}\,, 
  \label{primastab}
\Eeq
where we recall that $\Xi:=\calN(\dt\Rho)$.
Notice that \eqref{primastab} implies a part of~\eqref{contdepbis}.

\step
Second auxiliary estimate

We write the equation \eqref{prima} for both solutions 
and test the difference \aet\ by~$\Mu$.
The same we do with \eqref{seconda}, and test the difference by~$-\Mu$.
Then, we sum up and have, \aet,
\Bsist
  && \normaV\mu^2 + \normaVG\muG^2
  \non
  \\
  && = (\tauO - 1) \iO \dt\rho \, \mu
  + (\tauG - 1) \iG \dt\rhoG \, \muG
  + \iO \bigl( \rho_1 u_1 - \rho_2 u_2 \bigr) \cdot \nabla\mu
  \non
  \\
  && \quad {}
  + \iO \nabla\rho \cdot \nabla\mu
  + \iG \nablaG\rhoG \cdot \nablaG\muG
  \non
  \\
  && \quad {}
  + \iO \bigl( f'(\rho_1) - f'(\rho_2) \bigr) \mu
  + \iG \bigl( \fG'(\rho_{1\Gamma}) - \fG'(\rho_{2\Gamma}) \bigr) \muG \,.
  \non
\Esist
Now, we rearrange the \rhs\ and use the boundedness and the \Lip\ continuity of $f'$ and~$\fG'$,
as well as the \Holder\ and Young inequalities.
We obtain \aet~that
\Bsist
  && \normaV\mu^2 + \normaVG\muG^2
  \non
  \\
  && \leq \delta \, \normaH\mu^2
  + c_\delta \normaH{\dt\rho}^2
  + \delta \, \normaHG\muG^2
  + c_\delta \normaHG{\dt\rhoG}^2
  + \bigl(
    \norma\rho_6 \, \pier{{}\norma{u_1}_3 
    + \norma{\rho_2}_6{}} \, \norma u_3 
  \bigr) \norma{\nabla\mu}_2
  \non
  \\
  && \quad {}
  + \delta \, \normaH{\nabla\mu}^2
  + c_\delta \normaH{\nabla\rho}^2
  + \delta \, \normaHH{\nablaG\muG}^2
  + c_\delta \normaHG{\nablaG\rhoG}^2  
  \non
  \\
  && \quad {}
  + \delta \, \normaH\mu^2
  + c_\delta \normaH\rho^2
  + \delta \, \normaHG\muG^2
  + c_\delta \normaHG\rhoG^2\,,
  \non
\Esist
where $\delta>0$ is arbitrary.
By choosing $\delta$ small enough, using the Sobolev inequality,
and recalling that $\pier{u_1}\in\L\infty{\Lx3}$ and $\pier{\rho_2}\in\L\infty V$,
we deduce~that
\Beq
  \normaV\mu^2 + \normaVG\muG^2
  \leq c \bigl(
    \normaH{\dt\rho}^2
    + \normaHG{\dt\rhoG}^2
    \juerg{+ \normaV\rho^2}
    + \normaVG\rhoG^2
    + \norma u_3^2 
  \bigr) \quad \aet .
  \non
\Eeq
At this point, by accounting for~\eqref{primastab}, we conclude that
\Beq
  \norma\Mu_{\L\infty\calV}
  \leq c \, \norma u_{\H1{\Lx3}} \,. 
  \label{secondastab}
\Eeq

\step
Proof of Theorem~\ref{ContDepbis}

We \juerg{recall that \eqref{hpu} holds true} for both $u_1$ and $u_2$ and rewrite
the transport terms in \eqref{prima} in the form 
$\iO \nabla\rho_i\cdot u_i\,v$.
Then we take the difference of the equations, written for both solutions,
and apply Lemma~\ref{Elliptic} \aat\ with $\gamma=0$
and the following choice of $g$ and~$\gG$:
\Beq
  g = \bigl( {-\dt\rho} - \nabla\rho_1\cdot u_1 + \nabla\rho_2 \cdot u_2 \bigr) (t)
  = \bigl( {-\dt\rho} - \nabla\rho_1\cdot u + \nabla\rho \cdot u_2 \bigr) (t)
  \aand
  \gG = -\dt\rhoG(t) \,.
  \non
\Eeq
We then obtain that
\Beq
  \normaWW{\Mu(t)}
  \leq c \bigl( \normaVV{\Mu(t)} + \norma{\dt\rho(t)}_2 + \norma{u(t)}_3 + \norma{\nabla\rho(t)}_6 + \norma{\dt\rhoG(t)}_2\,,
  \non
\Eeq
where $c$ depends only on $\Omega$ and the norms of \pier{$\nabla \rho_1$} and $u_2$
in the spaces $\L\infty{\Lx6}$ and $\L\infty{\Lx3}$, respectively.
By combining this with \accorpa{primastab}{secondastab}, we deduce~that
\Beq
  \norma\Mu_{\L\infty\calW}
  \leq c \, \norma u_{\H1{\Lx3}}\,,
  \non
\Eeq
which is a part of~\eqref{contdepbis}.
In order to prove the remaining part of the estimate,
we write \eqref{seconda} for both solutions, take the difference,
and apply Lemma~\ref{Elliptic} \aat\ with $\gamma=0$ and \juerg{the choice}
\Beq
  g = \bigl( {-\tauO} \dt\rho - f'(\rho_1) + f'(\rho_2) + \mu \bigr) (t)
  \aand
  \gG = \bigl( {-\tauG} \dt\rhoG \juerg{-f_\Gamma'(\rho_{1\Gamma})
  +f'_\Gamma(\rho_{2\Gamma})} + \muG \bigr) (t) .
  \non
\Eeq
\juerg{We then obtain that}
\Beq
  \norma\Rho_{\L\infty\calW}
  \leq c \bigl(
    \norma\Rho_{\L\infty\calV}
    + \|\juerg{(g,\gG)}\|_{\L\infty\calH}
  \bigr)
  \leq c \, \norma u_{\H1{\Lx3}} \,,
  \non
\Eeq
\juerg{where the last inequality follows from}  \eqref{primastab} and \eqref{secondastab}.
\juerg{With this}, \eqref{contdepbis} is completely proved.\QED


\section*{Acknowledgments}
PC and GG gratefully acknowledge some financial support from 
the MIUR-PRIN Grant 2015PA5MP7 ``Calculus of Variations'' and 
the GNAMPA (Gruppo Nazionale per l'Analisi Matematica, 
la Probabilit\`a e le loro Applicazioni) of INdAM (Isti\-tuto 
Nazionale di Alta Matematica).


\vspace{3truemm}

{\small%

\Begin{thebibliography}{10}

\bibitem{bai}
F.~Bai, C.M. Elliott, A.~Gardiner, A.~Spence and A.M. Stuart,  
The viscous {C}ahn-{H}illiard equation. {I}. {C}omputations, 
{\it Nonlinearity\/}  {\bf 8} (1995) 131-160.

\bibitem{Barbu}
V. Barbu,
``Nonlinear Differential Equations of Monotone Type in Banach Spaces'',
Springer,
London, New York, 2010.

\bibitem{BL}
J. Bergh and J. L\"ofstr\"om,
``Interpolation spaces. An introduction'',
Grundlehren der mathematischen Wissenschaften {\bf 223},
Springer-Verlag, Berlin-New York, 1976.

\bibitem{Brezis}
H. Brezis,
``Op\'erateurs maximaux monotones et semi-groupes de contractions
dans les espaces de Hilbert'',
North-Holland Math. Stud.
{\bf 5},
North-Holland,
Amsterdam,
1973.

\bibitem{CahH} 
J.W. Cahn and J.E. Hilliard, 
Free energy of a nonuniform system I. Interfacial free energy, 
{\it J. Chem. Phys.\/}
{\bf 2} (1958) 258-267.

\bibitem{CaCo}
L. Calatroni and P. Colli,
Global solution to the Allen--Cahn equation with singular potentials and dynamic boundary conditions,
{\it Nonlinear Anal.\/} {\bf 79} (2013) 12-27.

\bibitem{CGM13}
 L.\ {C}herfils, S.\ {G}atti and A.\ {M}iranville, 
\newblock A variational approach to a {C}ahn--{H}illiard model in a domain 
with nonpermeable walls,
\newblock J.\ Math.\ Sci.\ (N.Y.) {\bf 189} (2013) 604-636. 

\bibitem{CMZ11}
 L.\ {C}herfils, A.\ {M}iranville and S.\ {Z}elik,
The {C}ahn--{H}illiard equation with logarithmic potentials,
{\it Milan J. Math.\/} {\bf 79} (2011) 561-596.%

\bibitem{CP14}
L.\ {C}herfils and M.\ {P}etcu, 
	\newblock A numerical analysis of the {C}ahn--{H}illiard equation with non-permeable walls,
	\newblock Numer.\ Math. {\bf 128} (2014) 518-549. 
	
\bibitem{CFP} 
R. Chill, E. Fa\v sangov\'a and J. Pr\"uss,
Convergence to steady states of solutions of the Cahn-Hilliard equation with dynamic boundary conditions,
{\it Math. Nachr.\/} 
{\bf 279} (2006) 1448-1462.

\bibitem{CFGS1}
P. Colli, M.H. Farshbaf-Shaker, G. Gilardi and J. Sprekels,
Optimal boundary control of a viscous Cahn--Hilliard system 
with dynamic boundary condition and double obstacle potentials, 
\pier{{\it SIAM J. Control Optim.\/} {\bf 53} (2015) 2696-2721.}

\bibitem{CFGS2}
P. Colli, M.H. Farshbaf-Shaker, G. Gilardi and J. Sprekels,
Second-order analysis of a boundary control problem 
for the viscous Cahn--Hilliard equation with dynamic boundary 
condition, {\it Ann. Acad. Rom. Sci. Ser. Math. Appl.\/} {\bf 7} 
(2015) 41-66.

\bibitem{CFS}
P. Colli, M.H. Farshbaf-Shaker and J. Sprekels,
A deep quench approach to the optimal control of an Allen--Cahn equation
with dynamic boundary conditions and double obstacles, 
{\it Appl. Math. Optim.\/} {\bf 71} (2015) 1-24.

\bibitem{CF1} 
P.\ {C}olli and T.\ {F}ukao, 
{C}ahn--{H}illiard equation with dynamic boundary conditions 
and mass constraint on the boundary, 
{\it J. Math. Anal. Appl.\/} {\bf 429} (2015) 1190-1213.

\bibitem{CF2} 
P.\ {C}olli and T.\ {F}ukao, 
Equation and dynamic boundary condition of 
Cahn--Hilliard type with singular potentials, 
{\it Nonlinear Anal.\/} {\bf 127} (2015) 413-433.

\bibitem{CGPS3} 
P. Colli, G. Gilardi, P. Podio-Guidugli and J. Sprekels,
Well-posedness and long-time behaviour for 
a nonstandard viscous Cahn-Hilliard system, 
{\it SIAM J. Appl. Math.} {\bf 71} (2011) 1849-1870.

\bibitem{CGS3}
P. Colli, G. Gilardi and J. Sprekels,
On the Cahn--Hilliard equation with dynamic 
boundary conditions and a dominating boundary potential,
{\it J. Math. Anal. Appl.\/} {\bf 419} (2014) 972-994.

\bibitem{CGS5}
P. Colli, G. Gilardi and J. Sprekels,
A boundary control problem for the pure Cahn-Hilliard equation
with dynamic boundary conditions,
\pier{{\it Adv. Nonlinear Anal.\/}} {\bf 4} (2015) 311-325. 

\bibitem{CGS4}
P. Colli, G. Gilardi and J. Sprekels,
A boundary control problem for the viscous Cahn-Hilliard equation
with dynamic boundary conditions,
{\it Appl. Math. Optim.\/} \pier{{\bf 73} (2016) 195-225.}

\bibitem{CS}
P. Colli and J. Sprekels,
Optimal control of an Allen--Cahn equation 
with singular potentials and dynamic boundary condition,
{\it SIAM J. Control Optim.\/} {\bf 53} (2015) 213-234.%

\bibitem{EllSt} 
C.M. Elliott and A.M. Stuart, 
Viscous Cahn--Hilliard equation. II. Analysis, 
{\it J. Differential Equations\/} 
{\bf 128} (1996) 387-414.

\bibitem{EllSh} 
C.M. Elliott and S. Zheng, 
On the Cahn--Hilliard equation, 
{\it Arch. Rational Mech. Anal.\/} 
{\bf 96} (1986) 339-357.

\bibitem{FG} 
E. Fried and M.E. Gurtin, 
Continuum theory of thermally induced phase transitions based on an order 
parameter, {\it Phys. D} {\bf 68} (1993) 326-343.

\bibitem{FY} 
T. Fukao and N. Yamazaki, 
A boundary control problem for the equation 
and dynamic boundary condition of Cahn--Hilliard type, 
to appear in ``Solvability, 
Regularity, Optimal Control of Boundary Value Problems for PDEs'', 
P.~Colli, A.~Favini, E.~Rocca, G.~Schimperna, J.~Sprekels~(ed.), 
Springer INdAM Series.

\bibitem{GiMiSchi} 
G. Gilardi, A. Miranville and G. Schimperna,
On the Cahn--Hilliard equation with irregular potentials and dynamic boundary conditions,
{\it Commun. Pure Appl. Anal.\/} 
{\bf 8} (2009) 881-912.

\bibitem{GiMiSchi2} 
G. Gilardi, A. Miranville and G. Schimperna,
Long-time behavior of the Cahn-Hilliard 
equation with irregular potentials 
and dynamic boundary conditions,
{\it Chin. Ann. Math. Ser. B\/} 
{\bf 31} (2010) 679-712.

\bibitem{GMS11} 
	 G.R.\ {G}oldstein, A.\ {M}iranville and G.\ {S}chimperna, 
	\newblock A {C}ahn--{H}illiard model in a domain with non-permeable walls,
	\newblock Phys. D {\bf 240} (2011) 754-766.

\bibitem{GM13} 
	 G.R.\ {G}oldstein and A.\ {M}iranville,
	\newblock A {C}ahn--{H}illiard--{G}urtin model with dynamic boundary conditions,
	\newblock Discrete Contin.\ Dyn.\ Syst.\ Ser.\ S {\bf 6} (2013) 387-400.

\bibitem{Gu} 
M. Gurtin, 
Generalized Ginzburg-Landau and Cahn-Hilliard equations based on a microforce balance,
{\it Phys.~D\/} {\bf 92} (1996) 178-192.

\bibitem{HW1}
M. Hinterm\"uller and D. Wegner, Distributed optimal control of the 
Cahn--Hilliard system including the case of a double-obstacle 
homogeneous free energy density, {\it SIAM J.
Control Optim.} {\bf 50} (2012) 388-418.
\bibitem{HW2}
M. Hinterm\"uller and D. Wegner, Optimal control of a semi-discrete 
Cahn--Hilliard--Navier--Stokes system, {\it SIAM J. Control Optim.} 
{\bf 52} (2014) 747-772.

\bibitem{Kub12} 
	 M.\ {K}ubo,
	\newblock The {C}ahn--{H}illiard equation with time-dependent constraint, 
	\newblock Nonlinear Anal. {\bf 75} (2012) 5672-5685. 

\bibitem{Lions}
J.-L.~Lions,
``Quelques m\'ethodes de r\'esolution des probl\`emes
aux limites non lin\'eaires'',
Dunod; Gauthier-Villars, Paris, 1969.

\bibitem{MZ}
A. Miranville and S. Zelik,
Robust exponential attractors for {C}ahn-{H}illiard type
equations with singular potentials,
{\it Math. Methods Appl. Sci.} {\bf 27} (2004) 545-582.

\bibitem{NovCoh}
A.~Novick-Cohen, On the viscous {C}ahn-{H}illiard equation, in
``Material instabilities in continuum mechanics'' ({E}dinburgh, 1985--1986),
Oxford Sci. Publ., Oxford Univ. Press, New York, 1988, pp.~329-342.

\bibitem{Podio}
P. Podio-Guidugli, 
Models of phase segregation and diffusion of atomic species on a lattice,
{\it Ric. Mat.} {\bf 55} (2006) 105-118.

\bibitem{PRZ} 
J. Pr\"uss, R. Racke and S. Zheng, 
Maximal regularity and asymptotic behavior of solutions for the 
Cahn--Hilliard equation with dynamic boundary conditions,  
{\it Ann. Mat. Pura Appl.~(4)\/}
{\bf 185} (2006) 627-648.
\bibitem{RZ} 
R. Racke and S. Zheng, 
The Cahn--Hilliard equation with dynamic boundary conditions, 
{\it Adv. Differential Equations\/} 
{\bf 8} (2003) 83-110.
\bibitem{RoSp}
E. Rocca and J. Sprekels,
\pier{Optimal distributed control of 
a nonlocal convective Cahn--Hilliard equation by 
the velocity in three dimensions,
{\it SIAM J. Control Optim.\/} {\bf 53} (2015) 1654-1680.}

\bibitem{Simon}
J. Simon,
Compact sets in the space $L^p(0,T; B)$,
{\it Ann. Mat. Pura Appl.~(4)\/} 
{\bf 146}, (1987), 65-96.

\bibitem{Wang}
Q.-F. Wang, Optimal distributed control of nonlinear {C}ahn-{H}illiard 
systems with computational realization, 
{\it J. Math. Sci. (N. Y.)\/} {\bf 177} (2011) 440-458.
\bibitem{WaNa}
Q.-F. Wang and S.-I. Nakagiri, Optimal control of distributed parameter 
system given by Cahn--Hilliard equation, 
{\it Nonlinear Funct. Anal. Appl.} {\bf 19} (2014) 19-33.
\bibitem{WZ} H. Wu and S. Zheng,
Convergence to equilibrium for the Cahn--Hilliard equation 
with dynamic boundary conditions, {\it J. Differential Equations\/}
{\bf 204} (2004) 511-531.
\bibitem{ZL1}
X.P. Zhao and C. C. Liu, Optimal control of the convective Cahn--Hilliard equation, 
{\it Appl. Anal.\/} {\bf 92} (2013) 1028-1045.
\bibitem{ZL2}
X.P. Zhao and C. C. Liu, Optimal control for the convective Cahn--Hilliard equation 
in 2D case, {\it Appl. Math. Optim.} {\bf 70} (2014) 61-82.
\bibitem{ZW}
J. Zheng and Y. Wang, Optimal control problem for {C}ahn--{H}illiard equations
with state constraint, {\it J. Dyn. Control Syst.\/} {\bf 21} (2015) 257-272.%

\End{thebibliography}

}

\End{document}
